\documentclass{article}

\usepackage{amsmath}
\usepackage{amsfonts}
\usepackage{amssymb}
\usepackage{amsthm}
\usepackage[T1]{fontenc}
\usepackage[latin1]{inputenc}
\usepackage[greek,francais,english]{babel}
\usepackage{graphics}
\usepackage{vmargin}
\setmarginsrb{1in}{1in}{1in}{1in}{0cm}{0cm}{0cm}{1.5cm}

\newcommand{\R}{\mathbb R}
\newcommand{\N}{\mathbb N}

\newcommand{\K}{\mathbb K}
\newcommand{\e}{\varepsilon}
\newcommand{\be}{\begin{equation}}
\newcommand{\ee}{\end{equation}}
\newcommand{\dv}{\mathrm{div}}

\newcommand{\p}{\partial}
\newcommand{\ue}{u^{\varepsilon}}

\newcommand{\fe}{f^{\varepsilon}}
\newcommand{\geps}{g^{\varepsilon}}

\newcommand{\me}{m^{\varepsilon}}
\newcommand{\sgn}{\mathrm{sgn}}
\newcommand{\Si}{\sum_{i=1}^N}
\newcommand{\supp}{\mathrm{Supp}}
\newcommand{\Id}{\mathrm{Id}}
\newcommand{\bu}{\bar{u}}

\newcommand{\ds}{\displaystyle}

\newcommand{\mean}[1]{\left\langle #1\right\rangle}
\newcommand{\ta}{\tilde{a}_0}
\newcommand{\tta}{\tilde{a}}
\newcommand{\fl}{f_{\lambda}}
\newcommand{\ul}{u_{\lambda}}

\newtheorem{prop}{Proposition}

\newtheorem{thm}{Theorem}
\newtheorem{defi}{Definition}
\newtheorem{corol}{Corollary}
\newtheorem{rem}{Remark}
\begin{document}
\title{Homogenization of nonlinear scalar conservation laws}
\author{Anne-Laure Dalibard}

\bibliographystyle{amsplain}
\maketitle

\begin{abstract}
We study the limit as $\e\to 0$ of the entropy solutions of the
equation $\p_t \ue +
\dv_x\left[A\left(\frac{x}{\e},\ue\right)\right] =0$. We prove
that the sequence $\ue$ two-scale converges towards a function
$u(t,x,y)$, and $u$ is the unique solution of a limit evolution
problem. The remarkable point is that the limit problem is not a
scalar conservation law, but rather a kinetic equation in which
the macroscopic and microscopic variables are mixed. We also prove
a strong convergence result in $L^1_{\text{loc}}$.
\end{abstract}


\section{Introduction}

This article is concerned with the asymptotic behavior of the
sequence $\ue\in\mathcal C([0,\infty), L^1_{\text{loc}}(\R^N))$, as the parameter $\e$ vanishes, where $\ue$ is the entropy solution of
the scalar conservation law
\begin{gather} \frac{\p \ue(t,x)}{\p t} +
\sum_{i=1}^N \frac{\p }{\p x_i}A_i\left(\frac{x}{\e},\ue(t,x)\right)=0\quad t\geq 0,\  x\in\R^N,\label{eqdepart}\\
\ue(t=0)=u_0\left(x,\frac{x}{\e}  \right).\label{CIdepart}
\end{gather}

The functions $A_i=A_i(y, v)$ ($y\in \R^N,\ v\in\R$) are assumed to
be $Y$-periodic, where $Y=\Pi_{i=1}^N (0,T_i)$ is the unit cell, and $u_0$ is also assumed to be periodic in its second variable.

Under regularity hypotheses on the flux, namely $A\in
W^{2,\infty}_{\text{per,loc}}(\R^{N+1})$, and when the initial
data $\ue(t=0)$ belongs to $L^{\infty}$, it is known that there exists
a unique entropy solution $\ue$ of the above system for all $\e>0$
given (see \cite{dafermos,kruzkhov1,kruzkhov2,Serre}). The study of the homogenization of such
hyperbolic scalar conservation laws has been investigated by
several authors, see for instance \cite{wetransport,SCLWE,ESerre}, and in the linear case \cite{HouXin,JabinTzavaras}. In dimension one, there is also an
equivalence with Hamilton-Jacobi equations which allows to use the
results of \cite{LPV}. In general, the results obtained by these
authors can be summarized as follows: there exists a function
$u^0=u^0(t,x,y)$ such that \be \ue -u^0\left(t,x,\frac x \e\right)
\to 0 \quad \text{in }L^1_{\text{loc}}((0,\infty)\times \R^N).
\label{strongCV}\ee The function $u^0(t,x,y)$ satisfies a microscopic equation,
called cell problem, and an evolution equation, which is a scalar conservation law in which the
coefficients depend on the microscopic variable $y$. In general,
there is no ``decoupling'' of the macroscopic variables $t,x$, and
the microscopic variable $y$: the average of $u^0$ with
respect to the variable $y$ is not the solution of an ``average''
conservation law.

To our knowledge, there are no results as soon as the dimension is
strictly greater than one when the flux does not satisfy a
structural condition of the type $A(y,\xi)=a(y)g(\xi)$. Here, we
investigate the behavior of the family $\ue$ for arbitrary fluxes.
We prove that \eqref{strongCV} still holds in some sense which
will be precised later on, and the function $u^0$ is a solution of
a microscopic cell problem. Precisely, we prove that even though there is no simple evolution equation satisfied by the function $u^0$ itself, the function
$$
f(t,x,y,\xi)=\mathbf 1_{\xi<u^0}
$$
is the unique solution of a linear transport equation, with a source term which is a Lagrange multiplier accounting for the constraints on $f$. This statement is reminiscent of the kinetic formulation for scalar conservation laws (see \cite{LPTCRAS,LPT,artBP}, the general presentation in \cite{BP}, and
\cite{kinform} for the heterogeneous case); this is not surprising since our method of proof relies on the kinetic formulation for equation \eqref{eqdepart}. However, in general, it is unclear whether $u^0$ is the solution of a scalar conservation law. Thus the kinetic formulation appears as the ``correct'' vision of the entropy solutions of \eqref{eqdepart}, at least as far as homogenization is concerned.

The rest of this introduction is devoted to the presentation of our main results. We begin with the description of the asymptotic problem, and then we state the convergence results in the general case.


\subsection{Description of the asymptotic evolution problem}

We first introduce the asymptotic evolution problem, for which we state an existence and uniqueness result; then we explain how this asymptotic problem can be understood formally.

In the following, we set, for $(y,\xi)\in \R^{N+1}$,
\begin{gather*}
a_i(y,\xi)=\frac{\p A_i}{\p\xi}(y,\xi),\quad 1\leq i \leq N,\\
a_{N+1}(y,\xi)=-\dv_y A(y,\xi).
\end{gather*}
We set $a(y,\xi)=(a_1(y,\xi),\cdots, a_{N+1}(y,\xi))\in \R^{N+1}$. Notice that $\dv_{y,\xi}a(y,\xi)=0$. These notations were introduced in \cite{kinform}.

Before giving the definition of the limit system, we recall the kinetic formulation for equation \eqref{eqdepart}, which was derived in \cite{kinform}. Indeed, we believe it may shed some light on the limit system. Let $\ue$ be an entropy solution of \eqref{eqdepart}. Then there exists a non-negative measure $\me\in M^1 ((0,\infty)\times \R^{N+1})$ such that $\fe=\mathbf 1_{\xi<\ue(t,x)}$ is a solution of the transport equation
\begin{eqnarray}
&&\p_t\fe + a_i \left( \frac{x}{\e},\xi \right) \p_{x_i} \fe + \frac{1}{\e}a_{N+1} \left( \frac{x}{\e},\xi \right) \p_{\xi} \fe=\p_{\xi}\me,\label{eq:depart_kf_gal}\\
&&\fe(t=0,x,\xi)=\mathbf 1_{\xi<u_0\left( x,\frac{x}{\e} \right)}.
\end{eqnarray}
In fact, this equation was derived in \cite{kinform} for the function $\geps(t,x,\xi)=\chi(\xi,\ue(t,x))$, where $\chi(\xi,u)=\mathbf 1_{0<\xi<u}-\mathbf 1_{u<\xi<0} $, for $u,\xi\in\R$, and under the additional assumption $a_{N+1}(y,0)=0$ for all $y\in\R^N$. However, it is easily proved, using the identity $\fe=\geps + \mathbf 1_{\xi<0}$, that $\fe$ satisfies \eqref{eq:depart_kf_gal}, even when $a_{N+1}(y,0)$ does not vanish.

We now define the limit system, which is reminiscent of equation \eqref{eq:depart_kf_gal} :

\begin{defi}

Let $f\in L^{\infty}([0,\infty)\times \R^N\times Y \times\R)$, $u_0\in L^{\infty}(\R^N\times Y)$. We say that $f$ is a {\it generalized kinetic solution of the limit problem, with initial data $\mathbf 1_{\xi<u_0}$,} if there exists a distribution $\mathcal M\in \mathcal D'_{\text{per}}([0,\infty)\times \R^N \times Y \times \R)$ such that $f$ and $\mathcal M$ satisfy the following properties:

\begin{enumerate}
 \item Compact support in $\xi$: there exists a constant $M>0$ such that
\begin{gather}
\supp \;\mathcal M \subset [0,\infty)\times \R^N \times Y \times [-M,M],\label{hyp:compact_support1}\\
f(t,x,y,\xi)=1\quad \text{if }\xi<-M,\label{hyp:compact_support2}\\
 f(t,x,y,\xi)=0\quad \text{if }\xi>M.\label{hyp:compact_support3}
\end{gather}

\item Microscopic equation for $f$: there exists a non-negative measure $m\in M^1((0,\infty)\times \R^N\times Y\times \R)$ such that $f$ is a solution in the sense of distributions of
\be
\dv_{y,\xi} (a(y,\xi) f(t,x,y,\xi))=\p_{\xi} m,\label{eq:cell_kin}\ee

and $\supp\: m\subset[0,\infty)\times \R^N\times Y \times [-M,M]$.

\item Evolution equation: the couple $(f,\mathcal M)$ is a solution in the sense of distributions of
\be
\left\lbrace
    \begin{array}{l}
        \ds\p_t f + \Si a_i(y,\xi) \p_{x_i}  f= \mathcal M,\\
        \ds f(t=0,x,y,\xi)=\mathbf 1_{\xi<u_0(x,y)}=:f_0(x,y,\xi);
\end{array}
\right.
\label{eq:lim_evol_eq}
\ee
In other words, for any test function $\phi\in\mathcal D_{\text{per}}([0,\infty)\times \R^N \times Y \times \R)$,
\begin{multline*}
\int_0^\infty \int_{\R^N\times Y\times\R}f(t,x,y,\xi)\left\lbrace \p_t \phi(t,x,y,\xi)+ \Si a_i(y,\xi) \p_{x_i}  \phi (t,x,y,\xi)\right\rbrace \:dt\:dx\:dy\:d\xi =\\
=-\left\langle \phi,\mathcal M\right\rangle_{\mathcal D,\mathcal D'} -  \int_{\R^N\times Y\times\R}\mathbf 1_{\xi<u_0(x,y)}\phi(t=0,x,y,\xi)\:dx\:dy\:d\xi.
\end{multline*}
\item Conditions on $f$: there exists a non-negative measure $\nu\in M^1_{\text{per}}([0,\infty)\times \R^N \times Y \times \R)$ such that
\begin{gather}
\p_{\xi} f = - \nu,\label{hyp:dfdxi}\\
0\leq  f(t,x,y,\xi)\leq 1\quad \text{almost everywhere.}
\end{gather}
And for all compact set $K\subset \R^N$,
\be
\frac{1}{\tau}\int_0^\tau\left|\left|f(s)-f_0 \right|  \right|_{L^2(K\times Y\times \R)}\:ds\underset{\tau\to 0}{\longrightarrow}0\label{hyp:cont_time0}.
\ee
\item Condition on $\mathcal M$: define the set
\begin{eqnarray*}
 \mathcal G:=\left\{\psi\in L^\infty_{\text{loc}}(Y\times \R)\right.,&& \p_{\xi} \psi\geq 0, \text{and }\exists \mu\in M^1_{\text{per}}(Y\times \R), \ \exists  C>0, \exists\alpha_-\in\R,\\
&&\dv_{y,\xi}(a\psi)=-\p_{\xi }\mu, \ \supp \:\mu\subset Y\times [-C,C],\ \mu\geq 0,\\
&&\left.\psi(y,\xi)=\alpha_-\text{ if } \xi<-C\right\}.
\end{eqnarray*}

Then for all $\varphi\in\mathcal D([0,\infty)\times \R^N)$ such that $\varphi\geq 0$, the function $\mathcal M \ast_{t,x} \varphi$ belongs to $\mathcal C ([0,\infty)\times \R^N, L^2(Y\times \R))$, and
\be
\forall (t,x)\in[0,\infty)\times \R^N,\quad \forall \psi\in\mathcal G,\quad\int_{Y\times \R } \left(\mathcal M \ast_{t,x} \varphi\right)(t,x,\cdot) \:\psi\leq 0
.
\label{hyp:M}
\ee

\end{enumerate}

\label{def:lim_ev_eq}
\end{defi}

We now state an existence and uniqueness result for solutions of the limit problem :
\begin{thm}
Let $A\in W^{2,\infty}_{\text{per,loc}}(Y\times \R)$.

\begin{enumerate}
 \item Existence: let $u_0\in L^1_{\text{loc}}(\R^N; \mathcal C_{\text{per}}(Y))\cap L^\infty(\R^N)$ such that there exists a
non-negative measure $m_0=m_0(x,y,\xi)$ such that
$f_0(x,y,\xi)=\mathbf 1_{\xi<u_0(x,y)}$ is a solution of \be \Si
\frac{\p}{\p y_i}\left(a_i(y,\xi) f_0 \right) +\frac{\p}{\p
\xi}\left(a_{N+1}(y,\xi) f_0 \right)=\frac{\p m_0}{\p
\xi}\label{eq:kinform_stat} \ee
and $\supp\:m_0 \subset \R^N\times Y\times [-M,M]$, where $M=||u_0||_{\infty}$.

Assume that there exists $u_1$, $u_2\in L^{\infty}(Y)$ such that
$\mathbf 1_{\xi< u_i}$ is a solution of \eqref{eq:kinform_stat}
for $i=1,2$, for some non-negative measures $m_1$, $m_2$, and \be
u_1(y)\leq u_0(x,y)\leq u_2(y)\quad \text{for a.e. }x\in\R^N,y\in
Y. \label{hyp:u0u1u2}\ee

Then there exists a generalized kinetic solution $f$ of the limit problem (in the sense of definition \ref{def:lim_ev_eq}), with initial data $f_0$.

\item ``Rigidity'': let $u_0\in L^\infty(\R^N\times Y)$, and let $f\in L^{\infty} ([0,\infty)\times\R^N\times Y \times \R)$ be a generalized kinetic solution of the limit problem, with initial data $f_0=\mathbf 1_{\xi<u_0}$. Then there exists a function $u\in L^\infty([0,\infty)\times\R^N\times Y)$ such that
$$
f(t,x,y,\xi)=\mathbf 1_{\xi<u(t,x,y)}\quad \text{almost everywhere.}
$$

\item Uniqueness and contraction principle: let $u_0, v_0\in L^\infty(\R^N\times Y)$, and let $f,g$ be two generalized kinetic solutions of the limit problem with initial data $\mathbf 1_{\xi<u_0}$ and $\mathbf 1_{\xi<v_0}$ respectively.
Then there exists a constant $C>0$ such that for all $t>0$, for all $R, R'>0$,
\be
||f(t) -g(t)||_{L^1(B_R\times Y\times\R)}\leq e^{Ct+R}\left(||u_0 -v_0||_{L^1(B_{R'}\times Y)} + e^{-R'}\right).\label{in:contraction}
\ee

\end{enumerate}
As a consequence, for all $u_0\in L^\infty(\R^N\times Y)\cap L^1_{\text{loc}}(\R^N, \mathcal C_{\text{per}} (Y))$ satisfying \eqref{eq:kinform_stat} and \eqref{hyp:u0u1u2}, there exists a unique generalized kinetic solution $f\in L^{\infty} ([0,\infty)\times\R^N\times Y \times \R)$ of the limit problem.
\label{thm:lim_ev_eq}
\end{thm}

\begin{rem}
Notice that for any function $v\in L^\infty(Y)$, $v$ is an entropy solution of the cell problem
$$
\dv_y A(y,v(y))=0
$$
if and only if there exists a non-negative measure $m\in M^1_{\text{per}}(Y\times \R)$ such that the equation
$$
\dv_{y,\xi}(a(y,\xi)\mathbf 1_{\xi<v(y)})=\p_{\xi} m.
$$
is satisfied in the sense of distributions on $Y\times \R$.

In the case where $A$ is divergence-free this condition becomes
$$
\Si \p_{y_i}(a_i(y,\xi)\mathbf 1_{\xi<v(y)})=0.
$$

Indeed, in that case, $v$ satisfies
$$
\Si \p_{y_i}(a_i(y,\xi)\mathbf 1_{\xi<v(y)})=\p_{\xi} m
$$
for some non-negative measure $m$ such that $\supp\: m \subset Y\times (-M,M)$. Consequently,
$$
\Si \p_{y_i}\left(\int_{-M- 1}^\xi a_i(y,w)\mathbf 1_{w<v(y)}\:dw \right)= m(y,\xi)\geq 0.
$$
Since the left-hand side has zero mean-value on $Y$ for all $\xi\in[-M,M]$, we deduce that $m=0$. Thus, in the case where the flux $A$ is divergence free, the limit system takes a slightly simpler form: conditions \eqref{eq:cell_kin}, \eqref{hyp:M} become
\begin{gather}
\ds \dv_y (a(y,\xi) f(t,x,y,\xi))=0,\nonumber\\
\p_t f + \Si a_i(y,\xi)\p_{x_i} f = \mathcal M,\nonumber\\
\left\lbrace
\begin{array}{l}
\int_{Y\times \R } \left(\mathcal M \ast_{t,x} \varphi\right)(t,x,\cdot) \:\psi\leq 0,\\
\forall \psi\in L^{\infty}_{\text{loc}}(Y\times \R),\ \dv_y(a\psi)=0,\text{ and }\p_{\xi} \psi\geq 0.
\end{array}
\right.\label{hyp:Mdivfree}
\end{gather}
All the other properties remain the same.

\label{rem:lim_syst_divfree}

\end{rem}

\begin{rem}
Assume that the flux $A$ is divergence-free, and set
\begin{gather*}
C_1:=\{\psi\in L^2_{\text{loc}}(Y\times\R),\ \Si \frac{\p}{\p y_i}\left(a_i(y,\xi) \psi(y,\xi)  \right)=0\}, \\
C_2:=\{\psi\in L^\infty_{\text{loc}}(Y\times\R),\ \p_{\xi}\psi \geq 0\}.
\end{gather*}
Then $C_1$, $C_2$ are convex sets of the vector space $L^2_{\text{loc}}(Y\times\R)$; thus condition \eqref{hyp:M} can be re-written as : for all $\varphi\in\mathcal D((-\infty,0)\times \R^N)$ such that $\varphi\geq 0$, for all $(t,x)\in(0,\infty)\times \R^N$, we have
$$
\mathcal M\ast \varphi(t,x)\in (C_1\cap C_2)^{\circ},
$$
where $C^{\circ}$ denotes the normal cone of $C$. Let us recall that when the space dimension is finite (that is, if $C_1,C_2$ are convex cones in $\R^d$ for some $d\in\N$), then
$$
(\mathrm{cl}(C_1)\cap \mathrm{cl}(C_2))^{\circ}=\mathrm{cl}\left(C_1^{\circ} + C_2^{\circ}\right),
$$
where $\mathrm{cl}(A)$ denotes the closure of the set $A$.

If we forget about the closure and the fact that we are considering convex sets in an infinite dimensional space, then we are tempted to write
$$
\mathcal M\ast \varphi(t,x)\in (C_1\cap C_2)^{\circ}=\mu_1 + \mu_2,
$$
with $\mu_i\in C_i^{\circ}$, $i=1,2$. Moreover, very formally, we have
$$
C_2^{\circ}\text{``$=$''}\{\p_{\xi }m,\quad m\text{ non-negative measure}\}.
$$
Thus, we may think of $\mathcal M$ as some distribution of the form
$$
\mathcal M = \p_{\xi} m + \mu_1,
$$
with $m$ a non-negative measure on $[0,\infty)\times \R^N\times Y\times \R$, and $\mu_1\in C_1^{\circ}$.

Of course, these computations are not rigorous, but we believe they may help the reader understanding the action of the distribution $\mathcal M$ (at least in the divergence-free case), even though the precise structure of $\mathcal M$ shall not be needed in the proof. Inequality \eqref{hyp:M} is sufficient for all the applications in this paper.

\end{rem}

Let us stress that uniqueness for the limit problem holds, even though the cell problem does not have a unique solution in general; indeed, in the linear divergence free case, that is, if $A(y,\xi)=a(y)\xi$, with $\dv_y a=0$, then a function $u$ is a solution of the cell problem if
$$
\dv_y (a(y)u(y))=0,\quad \mean{u}_Y=0.
$$
The constant function equal to zero is a solution of this equation, but in general there are other entropy solutions: think for instance of the case where $N=2$, and
$$
a(y_1,y_2)=(-\p_2 \phi(y_1,y_2), \p_1 \phi(y_1,y_2)),
$$
for some function $\phi\in\mathcal C^2_{\text{per}}(Y)$. Then any function $u$ of the form $g(\phi) - \mean{g(\phi)}$, with $g$ a continuous function, is an entropy solution. Let us emphasize that nonlinearity assumptions on the flux are not enough to ensure uniqueness of solutions either, see for instance \cite{LPV}.

In Theorem \ref{thm:lim_ev_eq}, the uniqueness of the solution of the limit system derives from a contraction principle associated with the macroscopic evolution equation, rather than the microscopic cell equation. The well-preparedness of the initial data, that is, the fact that $u_0(x,\cdot)$ is an entropy solution of the cell problem, is fundamental.

On the other hand, the lack of uniqueness of solutions of the cell problem entails that in general, there is no notion of homogenized problem. Indeed, if $u$ is a solution of
$$
\dv_yA(y,p + u(y))=0,\quad \mean{u}_Y=0,
$$
\noindent then in general, the quantity
$$
\mean{A(\cdot, p+ u(\cdot) }
$$
depends on $u$ (except when $N=1$, and in some special cases, when $N=2$; see \cite{EvansGomes2,LPV}). Hence the macroscopic and microscopic scales cannot be decoupled: if $\mathbf 1_{\xi<u(t,x,y)}$ is a solution of the limit evolution problem, then $\bu(t,x)=\mean{u(t,x,\cdot)}$ does not satisfy any remarkable equation. This is the main consequence of the absence of uniqueness for the cell problem.

\vskip2mm

Let us mention an important particular case of the theorem \ref{thm:lim_ev_eq}, which we call the ``separate case''. We now assume that the
flux $A$ can be written $A(y,\xi)=a_0(y) g(\xi)$, with $\dv_y
a_0=0$. This case has already been thorougly investigated by
Weinan E in \cite{wetransport} in the case where $g'(\xi)\neq 0$
for all $\xi$, that is, when the function $g$ is strictly monotonous. Here, we
prove that his results hold with no restriction on $g$.

Let us introduce the so-called ``constraint space''
$$
\K_0:=\{ f\in L^1(Y); \ \dv_y(a_0f)=0\quad \text{in}\ \mathcal D'\},
$$
and the orthogonal projection $P_0$ on $\K_0\cap L^2(Y)$ for the
scalar product in $L^2(Y)$.

Then the following properties hold: for all $f,g\in L^2(Y)$, if
$f\in\K_0$, then
$$
P_0(fg)=fP_0(g).
$$
And if $f,g\in \K_0$, then the product $fg$ belongs to $\K_0$.
Notice also that all functions which do not depend on $y$ belong
to $\K_0$.

\begin{prop}
Assume that $u_0\in L^1(\R^N,\mathcal C_{\text{per}}(Y))\cap L^\infty(\R^N\times Y)$, and $u_0$
is such that $u_0(x,\cdot)\in \K_0$ for a.e. $x\in\R^N$.

Let $\ta=P_0(a_0)\in L^{\infty}(Y)$. Let $u=u(t,x;y)$ be the entropy
solution of the scalar conservation law \be \left\{
\begin{array}{l} \p_t u(t,x;y) + \dv_x\left(\ta(y)
g(u(t,x;y)) \right)=0,\quad t>0,\ x\in\R^N,\ y\in Y,\\
u(t=0,x;y)=u_0(x,y).
\end{array}\right.
\label{eq:homogenized_sep} \ee

Then the function $f(t,x,y,\xi)=\mathbf 1_{\xi<u(t,x,y)}$ is the unique generalized kinetic solution of the limit problem \eqref{eq:lim_evol_eq} with initial data $\mathbf 1_{\xi<u_0(x,y)}$. In that case the distribution $\mathcal M$ is given by
$$
\mathcal M = \frac{\p m}{\p \xi} + g'(\xi)(\ta(y)-a_0(y))\cdot \nabla_x f,
$$
where $m$ is the kinetic entropy defect measure associated with the function $u$, that is, $f$ is a solution of
$$
\p_t f + g'(\xi) \ta(y) \cdot \nabla_x f=\p_{\xi}m.
$$

As a consequence, the solution $u(t,x;y)$ of \eqref{eq:homogenized_sep} is an entropy solution of
$$
\dv_y A(y,u)=0
$$
for almost every $(t,x)\in (0,\infty)\times \R^N$.

\label{thm:homogenized_sep}
\end{prop}


\subsection{Convergence results}

Our first result is concerned with entropy solutions of \eqref{eqdepart}.

\begin{thm}
Let $A\in W^{2,\infty}_{\text{per,loc}}(\R^{N+1})$. Assume that the initial data $u_0\in L^1_{\text{loc}}(\R^N, \mathcal C_{\text{per}}(Y)$ satisfies \eqref{eq:kinform_stat}, \eqref{hyp:u0u1u2}. Let $f=\mathbf 1_{\xi<u}$ be the unique generalized kinetic solution of the limit problem, with initial data $\mathbf 1_{\xi<u}$; the existence of $f$ follows from theorem \ref{thm:lim_ev_eq}. Then as $\e$ vanishes,
\be
\mathbf 1_{\xi<\ue(t,x)} \stackrel{\text{2
sc.}}{\rightharpoonup}\mathbf 1_{\xi<u(t,x,y)}.\label{2sc_cv}
\ee

As a consequence, for all regularization kernels
$\varphi^{\delta}$ of the form
$$\varphi^{\delta}(x)=\frac 1{\delta^N} \varphi\left(\frac x{\delta} \right), \quad x\in\R^N,$$
with $\varphi\in\mathcal D(\R^N)$, $\int \varphi=1$, $0\leq
\varphi\leq 1$, we have, for all compact
$K\subset[0,\infty)\times\R^N$, \be \lim_{\delta\to 0} \lim_{\e\to
0}\left|\left |\ue(t,x) - u\ast_x\varphi^{\delta}\left(t,x,\frac
x\e \right)\right| \right|_{L^1(K)}=0. \label{strongCVregul}\ee

\label{thm:strg_cvg}

\end{thm}
\begin{rem}

The assumption \eqref{eq:kinform_stat} means that $u_0$ is
``well-prepared'' in the sense that $u_0(x,\cdot)$ is an entropy
solution of
$$
\dv_y\left(A(y,u_0(x,y)) \right)=0
$$
for a.e. $x\in\R^N$. If this hypothesis is not satisfied, then it
is expected that the behavior of the sequence $\ue$ will depend on
the nature of the flux. If the flux is linear, then oscillations will
propagate, and the cell equation \eqref{eq:cell_kin} shall not be
satisfied in general. If the flux satisfies some strong
nonlinearity assumption, on the contrary, the conjecture is that
the solution $\ue$ re-prepares itself in order to match the
microscopic profile dictated by the equation. Few results in this
direction are known in the hyperbolic case; the reader may consult for instance
\cite{AmadoriSerre,SCLWE,EngquistE,SerreHandbook}. In \cite{initiallayer}, the
author studies the same equation as \eqref{eqdepart} in which a
viscosity term of order $\e$ is added, and proves such a result,
but the method relies strongly on the parabolicity of the
equation.

\end{rem}

\begin{rem}

The way in which theorem \ref{thm:strg_cvg} is stated might seem slightly peculiar; indeed, convergence results of the type
$$
\ue- u\left(t,x,\frac x \e \right)\to 0\quad \text{in}\
L^1_{\text{loc}}
$$
are expected to hold. But in order to establish such a result, it seems necessary to prove that
$$
\lim_{\delta \to 0}\int_K\sup_{y\in Y}\left |u(t,x,y) - u\ast_x\varphi^{\delta}\left(t,x,y \right)\right|\:dt\:dx=0.$$
But the evolution equation for $u$ (or rather, for $\mathbf 1_{\xi<u}$) is given by definition \ref{def:lim_ev_eq}; since the distribution $\mathcal M$ allows for very few computations, it seems difficult to derive such estimates.
\end{rem}

The next result generalizes theorem \ref{thm:homogenized_sep} to weaker solutions of equation \eqref{eqdepart}, called kinetic solutions. In order to simplify the presentation, we explain how to generalize the result in the divergence-free case; we explain in the remark following the theorem how to derive an analogous result in the case where the flux $A$ is arbitrary.

Thus, for the reader's convenience, we first recall the definition of kinetic solutions in the divergence-free case (see \cite{kinform} for the heterogeneous case, and the presentation in \cite{BP} for the homogeneous case)~:

\begin{defi}[Kinetic solutions of \eqref{eqdepart}]
Let $\ue\in\mathcal C([0,\infty), L¹(\R^N))$.
Assume that there exists a non-negative measure $\me\in \mathcal C(\R_{\xi}, w-M¹([0,\infty)\times\R^N)$ such that for all $T>0$, the function
$$
\xi\mapsto\int_0^T\int_{\R^N}\me(t,x,\xi)\:dt\:dx
$$
is bounded on $\R$, and vanishes as $|\xi|\to\infty$.

\noindent Assume also that $\fe(t,x,\xi):=\chi(\xi,\ue(t,x))$ is a solution in the sense of distributions of the linear transport equation
\begin{gather} \frac{\p \fe}{\p t} +\Si
a_i\left(\frac{x}{\e},\xi\right)\p_{x_i} \fe=\frac{\p \me}{\p \xi}\quad t\geq 0,\  x\in\R^N,\label{eqdepartkf_div0}\\
\fe(t=0)=\chi\left(\xi,
u_0\left(x,\frac{x}{\e}\right)\right),\label{CIdepartkf_div0}
\end{gather}

\noindent Then it is said that $\ue$ is a kinetic solution of equation \eqref{eqdepart}.

\end{defi}
The existence of such solutions is only known when the flux satisfies additionnal regularity assumptions. Assume that $a_i\in \mathcal C^1_{\text{per}}(Y\times \R)$ for $1\leq i\leq N$, and assume
that there exists a constant $C$ such that \be |a(y,\xi)|\leq C
\left( 1 + |\xi|\right)\quad \forall y\in Y\ \forall \xi\in\R.
\label{hyp:growtha}\ee

\noindent Under such hypotheses, it is proved in \cite{kinform} that
for all $u_0\in L^1(\R^N, \mathcal C_{\text{per}}(Y))$, there
exists a unique function $\ue\in\mathcal C([0,\infty), L^1(\R^N))$
such that $\chi(\xi,\ue)$ is a solution of \eqref{eqdepart}; $\ue$
is called the kinetic solution of
\eqref{eqdepart}-\eqref{CIdepart}. And if $\ue$ is bounded in
$L^{\infty}((0,T)\times\R^N)$ for all $T>0$, then $\ue$ is the
entropy solution of \eqref{eqdepart}. Moreover, a contraction principle holds between kinetic solutions.

Let us now state the convergence result for kinetic solutions :

\begin{thm}
 Let $A\in W^{2,\infty}_{\text{per,loc}}(Y\times \R)$ such that $\dv_{y} A(y,\xi)=0$ for all $y,\xi$. Assume that $a_i\in \mathcal C^1_{\text{per}}(Y\times \R)$ for $1\leq i\leq N$, and that \eqref{hyp:growtha} is satisfied. Assume that the initial data $u_0$ belongs to $L^1(\R^N, \mathcal C^1_{\text{per}}(Y))$ and satisfies
$$
\Si\frac{\p}{\p y_i} \left( a_i(y,\xi) \chi(\xi,u_0)\right)=0.
$$

Let $\ue\in\mathcal C([0,\infty), L^1(\R^N))$ be the kinetic solution of \eqref{eqdepart} with initial data $u_0(x,x/\e)$. Then there exists a function $u\in L^\infty([0,\infty), L^1(\R^N\times Y))$ such that the convergence results \eqref{2sc_cv} and \eqref{strongCVregul} hold, and
$$
\frac{\p}{\p y_i}(a(y,\xi) \chi(\xi,u(t,x,y)))=0\quad \text{in }\mathcal D'.
$$

Moreover, if we set
$$
\mathcal M := \frac{\p}{\p t}\chi(\xi,u) + \Si a_i(y,\xi)\frac{\p}{\p x_i} \chi(\xi,u)\in\mathcal D',
$$
then $\mathcal M$ satisfies \eqref{hyp:Mdivfree}.

\label{thm:stcvkinsol}
\end{thm}

\begin{rem}

Let us explain how this result can be generalized to the case where the flux $A$ is arbitrary. First, the $L^1$ setting is not adapted to this case, because the $L^1$ norm is not conserved by the equation in general. Hence another notion of kinetic solutions is needed; the correct functional space should be of the type $V + L¹(\R^N)$, where $V$ is a fixed solution of the cell problem.

Then, the crucial point in Theorem \ref{thm:stcvkinsol} is to find a sequence $u_0^n$ such that $u_0^n$ converges towards $u_0$ in $L^1(\R^N, \mathcal C_{\text{per}}(Y))$, and for all $n\in\N$, $u_0^n $ satisfies \eqref{eq:kinform_stat}, \eqref{hyp:u0u1u2}. Finding such a sequence is easy in the divergence-free case, but seems more difficult in the general case, since solutions of the cell problem are not known. This seems to be the main obstacle to the generalization of Theorem \ref{thm:stcvkinsol} to arbitrary fluxes. If this step is admitted, it is likely that the proof of Theorem \ref{thm:stcvkinsol} can be adapted to general settings.

\end{rem}

The plan of the paper is the following:  in section
\ref{sec:pfstcv1} we prove, under the hypotheses of theorem \ref{thm:strg_cvg}, that the two-scale limit of the sequence $\mathbf 1_{\xi<\ue(t,x)}$ is a generalized kinetic solution of the limit system. In section \ref{sec:contraction}, we study the limit problem introduced in definition \ref{def:lim_ev_eq} and we prove the rigidity and uniqueness results in theorem \ref{thm:lim_ev_eq}; hence theorem \ref{thm:lim_ev_eq} and \ref{thm:strg_cvg} will be proved by the end of section \ref{sec:contraction}. In section \ref{sec:BGK}, we study a relaxation model of BGK type, approaching the limit system in the divergence free case. In section \ref{sec:separated}, we prove Proposition \ref{thm:homogenized_sep}.
Eventually, in section \ref{sec:furtherremarks}, we have gathered further remarks on the notion of limit evolution problem.

\section{Asymptotic behavior of the sequence $\ue$\label{sec:pfstcv1}}

In this section, we prove that the two-scale limit
of the sequence $\fe=\mathbf 1_{\xi<\ue(t,x)}$, say $f^0(t,x,y,\xi)$, is a generalized kinetic solution of the limit system; thus the existence result of Theorem \ref{thm:lim_ev_eq} follows from this section.
The organization is the following: we first derive
some basic (microscopic) properties for the function $f^0$. Then we explain how regularization by
convolution can be used in two-scale problems. The two other subsections are devoted to the other properties of the limit system, namely condition \eqref{hyp:M} and the strong continuity at time $t=0$.

\subsection{Basic properties of $f^0$}

We use the concept of two-scale convergence, formalized by G.
Allaire after an idea of G. N'Guetseng (see \cite{Allaire,NG}). The fundamental result in \cite{Allaire} can be
generalized to the present setting as follows:
\begin{corol}

Let $(\geps)_{\e>0}$ be a bounded sequence in
$L^{\infty}((0,\infty)\times\R^{N+1})$. Then there exists a
function $g^0\in L^{\infty}((0,\infty)\times\R^{N}\times Y \times
\R)$, and a subsequence $(\e_n)$ such that $\e_n\to 0$ as
$n\to\infty$, such that
$$
\int_0^\infty \int_{\R^{N+1}}g^{\e_n}(t,x,\xi) \psi\left(t,x,\frac
x \e,\xi \right)\:dt\:dx\:d\xi \to\int_0^\infty \int_{\R^{N}\times
Y \times \R}g^0(t,x,y,\xi)\psi(t,x,y,\xi)\:dt\:dx\:dy\:d\xi
$$
for all functions $\psi\in L^1((0,\infty)\times\R^{N+1}; \mathcal
C_{\text{per}}(Y))$.

It is said that the sequence $(g^{\e_n})_{n\in\N}$ two-scale converges towards $g^0$.

\label{cor:twosc}
\end{corol}

Here, the sequence $\fe$ is bounded by $1$ in $L^{\infty}$; hence
we can extract a subsequence, still denoted by $\e$, and find a
function $f^0\in L^{\infty}((0,\infty)\times\R^{N}\times Y \times
\R)$ such that $(\fe)$ two-scale converges to $f^0$. It is easily
checked that $f^0$ inherits the following properties from the
sequence $\fe$
\begin{gather}
0\leq f^0(t,x,y,\xi)\leq 1,\label{eq:sgnf0}\\
\p_{\xi} f^0=-\nu(t,x,y,\xi),\quad \nu\
\text{non-negative measure}.\label{eq:dxif0}
\end{gather}

Now, let us prove \eqref{hyp:compact_support2}-\eqref{hyp:compact_support3}: let
$$
M:= \max \left(||u_1||_{\infty} , ||u_2||_{\infty}\right),
$$
where $u_1,u_2$ are the functions appearing in assumption \eqref{hyp:u0u1u2}. Since $u_i\left( x/\e \right)$ is a stationary solution of \eqref{eqdepart}, by a comparison principle for equation \eqref{eqdepart}, we deduce that
$$
 u_1\left( \frac{x}{\e} \right)\leq \ue(t,x)\leq u_2\left( \frac{x}{\e} \right)\quad \text{for almost every }t>0,\ x\in\R^N.
$$
Thus $||\ue||_{L^\infty([0,\infty)\times \R^N}\leq M$, and for almost every $t,x,\xi$, for all $\e>0$,
\begin{gather*}
\fe(t,x,\xi)=1\quad \text{if } \xi<-M,\\
\fe(t,x,\xi)=0\quad \text{if } \xi>M.
\end{gather*}
Passing to the two-scale limit, we infer \eqref{hyp:compact_support2} and \eqref{hyp:compact_support3}.

Now, we derive a microscopic equation for $f^0$. First, multiplying \eqref{eq:depart_kf_gal} by $S'(\xi)$, with $S'\in\mathcal D(\R)$, and integrating on $(0,T)\times B_R\times \R$, with
$T>0$, $R>0$, yields
\begin{multline*}
\int_{B_R}\left(S(\ue(T,x)) -S\left(u_0\left(x,\frac x \e
\right)\right) \right)\:dx +\int_0^T \int_{\R}\int_{\p B_R}
a\left(\frac{x}{\e},\xi \right)\cdot n_R(x) \fe
S'(\xi)\:d\sigma_R(x)\:d\xi\:dt -\\- \frac{1}{\e}\int_0^T
\int_{\R}\int_{B_R}a_{N+1}\left(\frac{x}{\e},\xi \right)\fe
S''(\xi)\:dx\:d\xi\:dt=-\int_0^T \int_{\R}\int_{B_R} \me(t,x,\xi)
S''(\xi)\:dx\:d\xi\:dt,
\end{multline*}
where $n_R(x)$ is the outward-pointing normal to $B_R$ at a given point $x\in\p B_R$, and $d\sigma_R(x)$ is the Lebesgue measure on $\p B_R$.

Hence we obtain the following bound on $\me$
$$
\e\me((0,T)\times B_R\times \R)\leq C_{T,R}
$$
for all $\e>0$, $R>0$, $T>0$, and $\supp \:\me \subset
(0,\infty)\times \R\times [-M,M]$.

Consequently, there exists a further subsequence, still denoted by
$\e$, and a non-negative measure $m^0=m^0(t,x,y,\xi)$ such that $\e
\me$ two-scale converges to $m^0$ (the concept of two-scale convergence can easily be generalized to measures; the arguments are the same as in \cite{Allaire}, the only difference lies in the functional spaces). Moreover, $\supp\; m^0\subset
(0,\infty)\times \R\times Y\times  [-M,M]$.

We now multiply \eqref{eq:depart_kf_gal} by test functions of the type
$\e\varphi\left(t,x, x/\e,\xi \right)$, with $\varphi\in \mathcal D_{\text{per}}([0,\infty)\times \R^N\times Y\times \R)$,
and we pass to the two-scale limit. We obtain, in the sense of
distributions on $(0,\infty)\times\R^N\times Y\times \R$ \be
\frac{\p}{\p y_i}\left(a_i(y,\xi) f^0 \right) + \frac{\p}{\p
\xi}\left(a_{N+1}(y,\xi) f^0 \right)=\frac{\p m^0}{\p \xi}.
\label{eq:micro}\ee

\noindent Thus \eqref{eq:cell_kin} is satisfied, which  completes the derivation of the basic properties of $f^0$.

Now, we define the distribution
$$
\mathcal M := \frac{\p f^0}{\p t} + \Si a_i(y,\xi) \frac{\p   f^0}{\p{x_i}}.
$$
The distribution $\mathcal M$ obviously satisfies \eqref{hyp:compact_support1}. The next step is to prove that $\mathcal M$ satisfies \eqref{hyp:M}; since regularizations by convolution are involved in condition \eqref{hyp:M}, we now describe the links between convolution and two-scale convergence.


\subsection{Regularization by convolution and two-scale convergence\label{ssec:convol}}

In this subsection, we wish
to make a few remarks concerning the links between convolution and
two-scale convergence. Indeed, it is a well-known fact that if a
sequence $(f_n)$ weakly converges in $L^2(\R^N)$ towards a
function $f$, then for all convolution kernels
$\varphi=\varphi(x)$, the sequence $(f_n\ast \varphi)$ two-scale
converges in $L^2$ towards $f\ast \varphi$. It would be convenient
to have an analogue property for two-scale limits. However, in
general, if a sequence $\fe=\fe(x)$ is bounded in $L^2(\R^N)$ and
two-scale converges towards a function $f=f(x,y)\in
L^2(\R^N\times Y)$, then $\fe\ast \varphi$ does not two-scale
converge towards $f\ast_x \varphi$. Indeed, if $\psi=\psi(x,y)\in
L^2(\R^N, \mathcal C_{\text{per}}(Y))$, then
\begin{eqnarray*}
&&\int_{\R^N}\fe\ast \varphi(x)
\psi\left(x,\frac{x}{\e}\right)\:dx\\
&=&\int_{\R^{2N}}\fe(x')\varphi(x-x')\psi\left(x,\frac{x}{\e}\right)\:dx\:dx'\\
&=&\int_{\R^{N}}dx'\:\fe(x')\left[\int_{\R^N}\varphi(x-x')\psi\left(x,\frac{x}{\e}\right)\:dx\right].
\end{eqnarray*}
In general, the quantity between brackets in the last integral
cannot be written as a function of $x'$ and $x'/\e$, and it seems
difficult to pass to the limit as $\e\to 0$.

In order to get round this difficulty, let us suggest the
following construction, which is reminiscent of the doubling of
variables in the papers of Kruzkhov, see
\cite{kruzkhov1,kruzkhov2}. With the same notations as above,
consider the test function
$(\psi\ast_x\check{\varphi})\left(x,\frac{x}{\e} \right)$, where
$\check{\varphi}(x):=\varphi(-x)\ \forall x\in\R^N$. Then by
definition of the two-scale convergence,
$$
\int_{\R^N}\fe(x)\left[\psi\ast_x
\check{\varphi}\right]\left(x,\frac{x}{\e} \right)\: dx
\to\int_{\R^N\times Y}f(x,y)\left[\psi\ast_x
\check{\varphi}\right]\left(x,y \right)\: dx\:dy
$$
And
\begin{gather*}
\int_{\R^N}\fe(x)\left[\psi\ast_x
\check{\varphi}\right]\left(x,\frac{x}{\e} \right)\:
dx=\int_{\R^{2N}}\fe(x')\varphi(x-x')
\psi\left(x,\frac{x'}{\e}\right)\:dx\:dx',\\
\int_{\R^N\times Y}f(x,y)\left[\psi\ast_x
\check{\varphi}\right]\left(x,y \right)\: dx\:dy=\int_{\R^N\times
Y}\left[f\ast_x\varphi\right](x,y)\psi(x,y)\: dx\:dy.
\end{gather*}

Consequently, as $\e\to 0$, \be \int_{\R^{2N}}\fe(x')\varphi(x-x')
\psi\left(x,\frac{x'}{\e}\right)\:dx\:dx'\to\int_{\R^N\times
Y}\left[f\ast_x\varphi\right](x,y)\psi(x,y)\: dx\:dy
\label{lim:convol}\ee for all $\varphi\in\mathcal D(\R^N)$, for
all $\psi\in L^2(\R^N,\mathcal C_{\text{per}}(Y))$.

In fact, different assumptions on the function $\psi$ can be chosen; the key point is that $\psi$ should be an admissible test function in the sense of Allaire (see \cite{Allaire}). In particular, if there exist $\psi_1\in \mathcal D(\R^N)$, $\psi_2\in L^\infty(Y)$ such that
$$
\psi(x,y)=\psi_1(x) \psi_2(y),
$$
then $\psi$ is an admissible test function, and the limit \eqref{lim:convol} holds.

\subsection{Proof of the condition on $\mathcal M$\label{ssec:pf_hypM}}

The goal of this subsection is to prove that with
$$
\mathcal M= \p_t f^0 + \Si a_i(y,\xi) \p_{i} f^0,
$$
condition \eqref{hyp:M} holds; hence, let $\varphi\in\mathcal D(\R\times \R^N)$, $\theta\in \mathcal D(\R\times \R^N)$, such that
\begin{gather*}
\varphi\geq 0, \ \theta\geq 0,\\
\varphi(t,x)=0\ \forall t\geq 0\ \forall x\in\R^N,\quad \theta(t,x)=0\ \forall t\leq 0\ \forall x\in\R^N;
\end{gather*}
the function $\varphi$ shall be used as a convolution kernel, and $\theta$ as a test function, which explains the above hypotheses on the supports of $\varphi$ and $\theta$.

\noindent Let $\psi\in\mathcal G$ arbitrary (the definition of the set $\mathcal G$ is given in definition \ref{def:lim_ev_eq}). We have to prove that the quantity\begin{multline*}
                      A:=\int_0^\infty\int_0^\infty\int_{\R^{2N} \times Y\times \R} f^0(s,z,y,\xi)\: \left\{\p_t \varphi(t-s,x-z) + \Si a_i(y,\xi) \p_{i} \varphi(t-s,x-z)\right\}\times\\
\times\psi(y,\xi) \theta (t,x)\;d\xi\:dy\:dx\:dz\:ds\:dt                     \end{multline*}
in non-positive.

Before going into the technicalities, let us explain formally why the property is true; let us forget about the convolution and the regularity issues, and take the test function
$$
\theta(t,x) \psi\left( \frac{x}{\e},\xi \right)
$$
in equation \eqref{eq:depart_kf_gal}.

Let $R>\max(M,C)$; recall that $M$ and $C$ are such that $\supp \:f^0\subset [0,\infty) \times \R^N\times Y\times [-M,M]$, and $\psi(y,\xi)= \alpha_-$ if $\xi<-C$. Integrating on $[0,\infty) \times \R^N\times [-R,R]$, we obtain

\begin{eqnarray*}
&&\int_0^\infty \int_{\R^{N}}\int_{-R}^R\fe(t,x,\xi) \left[ \p_t \theta(t,x) + a_i \left( \frac{x}{\e},\xi \right)\p_{x_i}\theta(t,x) \right]\:\psi\left( \frac{x}{\e},\xi \right)\:dx\:d\xi\:dt\\
&& -\frac{1}{\e} \int_0^\infty \int_{\R^{N}}\int_{-R}^R \fe(t,x,\xi) \frac{\p \mu}{\p\xi}\left( \frac{x}{\e},\xi  \right)\theta(t,x)\:dx\:d\xi\:dt \\&&+ \alpha_- \int_0^\infty\int_{\R^{N}}\frac{1}{\e}a_{N+1} \left( \frac{x}{\e}, -R \right)\theta (t,x)\:dt\:dx \\
&=&\int_0^\infty \int_{\R^{N}}\int_{-R}^R\me(s,z,\xi) \p_{\xi}  \psi \left( \frac{x}{\e},\xi \right)\:dz\:d\xi\:ds-\int_{\R^{N}}\int_{-R}^R \mathbf 1_{\xi<u_0\left(x,\frac{x}{\e}\right)} \:\theta(t=0,x)\psi\left( \frac{x}{\e},\xi \right)\:dx\:d\xi .
\end{eqnarray*}

Notice that
$$
\frac{1}{\e}a_{N+1} \left( \frac{x}{\e}, -R \right)=-\dv_x A\left( \frac{x}{\e}, -R \right),
$$
and thus
\begin{eqnarray*}
&&\int_0^\infty \int_{\R^{N}}\int_{-R}^R\fe(t,x,\xi) \left[ \p_t \theta(t,x) + a_i \left( \frac{x}{\e},\xi \right)\p_{x_i}\theta(t,x) \right]\:\psi\left( \frac{x}{\e},\xi \right)\:dx\:d\xi\:dt\\
&=&\int_0^\infty \int_{\R^{N}}\int_{-R}^R\left[\me(s,z,\xi) \p_{\xi} \psi \left( \frac{x}{\e},\xi \right) - \frac{1}{\e}\mu\left( \frac{x}{\e},\xi \right)  \p_{\xi}\fe(t,x,\xi)\right]\theta(t,x) \:dz\:d\xi\:ds\\
&&- \alpha_-\int_0^\infty \int_{\R^{N}}A_i\left( \frac{x}{\e},-R \right)\p_{i} \theta(t,x)\:dt\:dx-\int_{\R^{N}}\int_{-R}^R\mathbf 1_{\xi<u_0\left( x,\frac{x}{\e}\right)}\:\theta(t=0,x)\psi\left( \frac{x}{\e},\xi \right)\:dx\:d\xi \\
&\geq & - \alpha_-\int_0^\infty \int_{\R^{N}}A_i\left( \frac{x}{\e},-R \right)\p_{i} \theta(t,x)\:dt\:dx-\int_{\R^{N}}\int_{-R}^R\mathbf 1_{\xi<u_0\left( x,\frac{x}{\e}\right)} \:\theta(t=0,x)\psi\left( \frac{x}{\e},\xi \right)\:dx\:d\xi.
\end{eqnarray*}
Passing to the limit as $\e\to 0$, we retrieve
\begin{eqnarray*}
&&\int_0^\infty \int_{\R^{N}}\int_{-R}^Rf^0(t,x,y,\xi)\left[ \p_t \theta(t,x) + a_i \left( y,\xi \right)\p_{x_i}\theta(t,x) \right]\:\psi\left( y,\xi \right)\:dx\:dy\:d\xi\:dt\\
&\geq &- \alpha_-\int_0^\infty \int_{\R^{N}\times Y}A_i(y,-R)\p_{i} \theta(t,x)\:dt\:dx -\int_{\R^{N}}\int_{-R}^R\mathbf 1_{\xi<u_0( x,y)} \:\theta(t=0,x)\psi\left( y,\xi \right)\:dx\:d\xi\\
&=&-\int_{\R^{N}}\int_{-R}^R\mathbf 1_{\xi<u_0( x,y)} \:\theta(t=0,x)\psi\left( y,\xi \right)\:dx\:d\xi.
\end{eqnarray*}
This means exactly that
$$
\frac{\p}{\p t}\int_{Y\times \R} f^0 \psi + \frac{\p}{\p x_i} \int_{Y\times \R}a_i f^0\psi\leq 0,
$$
or in other words, that $\int_{Y\times \R} \mathcal M \psi\leq 0$ in the sense of distributions on $[0,\infty)\times \R^N$.

Now, we go back to the regularizations by convolution. According to the preceding subsection,
\begin{multline*}
A=\lim_{\e\to 0}\int_0^\infty\int_0^\infty\int_{\R^{2N} \times Y\times \R} \fe(s,z,\xi)\: \left\{\p_t \varphi(t-s,x-z) + \Si a_i\left(\frac{z}{\e},\xi\right) \p_{i} \varphi(t-s,x-z)\right\}\times\\
\times\psi\left(\frac{z}{\e},\xi\right) \theta (t,x)\;d\xi\:dx\:dz\:ds\:dt.
   \end{multline*}
Hence, in \eqref{eq:depart_kf_gal}, we consider the test function
$$
\phi(s,z,\xi)=\left[\int_0^\infty \int_{\R^N}\varphi(t-s,x-z)\:\theta(t,x)\: dt\:dx \right] \: \psi_{\delta}\left( \frac{z}{\e},\xi \right) K(\xi),
$$
where $K$ is a cut-off function such that $0\leq K\leq 1$, $K\in\mathcal D(\R)$, $K(\xi)=1$ if $|\xi|\leq R$, and
$$
\psi_{\delta}:=\psi\ast_y\varphi_1^{\delta}\ast_{\xi}\varphi_2^{\delta},
$$
with $\varphi_1\in\mathcal D(\R^N)$, $\varphi_2\in\mathcal D(\R)$, $0\leq \varphi_i\leq 1$,
$\int \varphi_i=1$ for $i=1,2$, and
$$
 \varphi_1^{\delta}(y)= \frac{1}{\delta^N} \varphi_1\left(\frac{y}{\delta}  \right),\quad
 \varphi_2^{\delta}(\xi)= \frac{1}{\delta} \varphi_2\left(\frac{\xi}{\delta}  \right).
$$

According to \eqref{eq:depart_kf_gal}, we have
\begin{eqnarray}
&&\int_0^\infty \int_{\R^{N+1}}\fe(s,z,\xi) \left[ \p_s \phi(s,z,\xi) + \Si a_i \left( \frac{z}{\e},\xi \right)\p_{z_i}\phi(s,z,\xi) \right]\:dz\:d\xi\:ds\nonumber\\
&+& \frac{1}{\e}\int_0^\infty \int_{\R^{N+1}}\fe(s,z,\xi) a_{N+1}\left( \frac{z}{\e},\xi \right)\p_{\xi}\phi(s,z,\xi)\:dz\:d\xi\:ds\label{eq:kf_distrib}\\
&-&\int_0^\infty \int_{\R^{N+1}}\me(s,z,\xi) \p_{\xi} \phi(s,z,\xi)\:dz\:d\xi\:ds+ \int_{\R^{N+1}}\chi\left(\xi,u_0\left( z,\frac{z}{\e} \right)\right)\:\phi(s=0,z,\xi)\:dz\:d\xi\nonumber\\&=&0.\nonumber
\end{eqnarray}

\noindent And
\begin{eqnarray*}
 \p_s \phi(s,z,\xi)&=&- \left[\int_0^\infty \int_{\R^N}\p_t\varphi(t-s,x-z)\:\theta(t,x)\: dt\:dx \right] \: \psi_{\delta}\left( \frac{z}{\e},\xi \right) K(\xi),\\
\nabla_z \phi(s,z,\xi)&=&-\left[\int_0^\infty \int_{\R^N}\nabla_x\varphi(t-s,x-z)\:\theta(t,x)\: dt\:dx \right] \: \psi_{\delta}\left( \frac{z}{\e},\xi \right) K(\xi)\\
&&+ \frac{1}{\e}\left[\int_0^\infty \int_{\R^N}\varphi(t-s,x-z)\:\theta(t,x)\: dt\:dx \right] \:\left(\nabla_y \psi_{\delta}\right)\left( \frac{z}{\e},\xi \right) K(\xi),\\
\p_\xi \phi(s,z,\xi)&=&\left[\int_0^\infty \int_{\R^N}\varphi(t-s,x-z)\:\theta(t,x)\: dt\:dx \right]  K(\xi)\: \p_{\xi}\psi_{\delta}\left( \frac{z}{\e},\xi \right)\\
&& +\left[\int_0^\infty \int_{\R^N}\varphi(t-s,x-z)\:\theta(t,x)\: dt\:dx \right]\psi_{\delta}\left( \frac{z}{\e},\xi \right)\p_{\xi} K(\xi)  \\
 \phi(s=0,z,\xi)&=&\left[\int_0^\infty \int_{\R^N}\p_t\varphi(t,x-z)\:\theta(t,x)\: dt\:dx \right] \: \psi_{\delta}\left( \frac{z}{\e},\xi \right) K(\xi)=0.
\end{eqnarray*}
Thanks to the assumption on the support of $\varphi$, and the fact that
$$
\p_{\xi }\psi_{\delta}= \left(\p_{\xi}\psi\right) \ast_y\varphi_1^{\delta}\ast_{\xi}\varphi_2^{\delta}\geq 0,
$$
we have
$$
\left[\int_0^\infty \int_{\R^N}\varphi(t-s,x-z)\:\theta(t,x)\: dt\:dx \right]  K(\xi)\: \p_{\xi}\psi_{\delta}\left( \frac{z}{\e},\xi \right)\geq 0.
$$
Moreover, thanks to \eqref{hyp:compact_support2}, \eqref{hyp:compact_support3}, and the assumptions on $\psi$ and $K$, we have $\p_{\xi} K =0$ on $\supp\: \me$, and
\begin{eqnarray*}
 &&\left[\int_0^\infty \int_{\R^N}\varphi(t-s,x-z)\:\theta(t,x)\: dt\:dx \right]\psi_{\delta}\left( \frac{z}{\e},\xi \right)\p_{\xi} K(\xi) \fe(s,z,\xi)\\
&=&\alpha_- \left[\int_0^\infty \int_{\R^N}\varphi(t-s,x-z)\:\theta(t,x)\: dt\:dx \right]\p_{\xi} K(\xi).
\end{eqnarray*}

Hence, we obtain, for all $\e,\delta>0$,
\begin{eqnarray*}
&&-\int \fe(s,z,\xi)\left\{\p_t \varphi(t-s,x-z) + \Si a_i\left(\frac{z}{\e},\xi\right) \p_{i} \varphi(t-s,x-z)\right\}\times\\
&&\qquad \qquad\qquad\qquad \qquad\qquad\qquad \qquad\qquad\qquad \qquad\qquad\times\psi_{\delta}\left(\frac{z}{\e},\xi\right) \theta (t,x)\:d\xi\:dx\:dz\:ds\:dt\\
&& + \frac{1}{\e}\int\fe(s,z,\xi) \:a \left( \frac{z}{\e},\xi \right)\cdot \nabla_{y,\xi} \psi_{\delta}\left( \frac{z}{\e},\xi \right) \varphi(t-s,x-z)\:\theta(t,x)\:K(\xi)\: dt\:dx \:ds\:dz\:d\xi\\
&&+\frac{\alpha_-}{\e}  \int\varphi(t-s,x-z)\:\theta(t,x)\:\p_{\xi} K(\xi) a_{N+1}\left( \frac{z}{\e},\xi \right) \: dt\:dx \:ds\:dz\:d\xi \\
&\geq &0.
\end{eqnarray*}

Following the formal calculations above, we have to investigate the sign of the term
$$
\int\fe(s,z,\xi)\: a \left( \frac{z}{\e},\xi \right)\cdot \nabla_{y,\xi} \psi_{\delta}\left( \frac{z}{\e},\xi \right) \varphi(t-s,x-z)\:\theta(t,x)\:K(\xi)\: dt\:dx \:ds\:dz\:d\xi.
$$
Since $\dv_{y,\xi}(a\psi)=-\p_{\xi}\mu$, we have
$$
\dv_{y,\xi}(a \psi_{\delta})=-\frac{\p \mu_{\delta}}{\p
\xi} + r_{\delta}
$$
where $\mu_{\delta}= \mu\ast_y\varphi_1^{\delta}\ast_{\xi}\varphi_2^{\delta}$. Then
\begin{eqnarray*}
 &&-\int_0^\infty \int_{\R^{N+1}}\fe(s,z,\xi) \frac{\p  \mu_{\delta}}{\p \xi }\left(\frac{x}{\e} ,\xi \right) \left[\int_0^\infty \int_{\R^N}\varphi(t-s,x-z)\:\theta(t,x)\: dt\:dx \right]\:ds\:dz\:d\xi\\
&=&-\int_0^\infty \int_{\R^{N+1}}\delta(\xi=\ue(t,x)) \mu_{\delta}\left(\frac{x}{\e} ,\xi \right) \left[\int_0^\infty \int_{\R^N}\varphi(t-s,x-z)\:\theta(t,x)\: dt\:dx \right]\:ds\:dz\:d\xi\leq 0.
\end{eqnarray*}
Hence, we have to prove that as $\delta\to 0$,
$$
r_{\delta}\to 0  \quad \text{in }L^1_{\text{loc}}(Y\times \R).
$$
The proof is quite classical. We have
\begin{eqnarray*}
r_{\delta}(y,\xi)&=& a(y,\xi)\psi\ast\left(\nabla_{y,\xi}\varphi_1^{\delta}\varphi_2^{\delta}\right)
-\left[a(y,\xi)\psi\right]\ast\left(\nabla_{y,\xi}\varphi_1^{\delta}\varphi_2^{\delta}\right)\\
&=&\Si\int\left[a_i(y,\xi)-a_i(y_1,\xi_1)\right]\psi(y_1,\xi_1)\p_{y_i}\varphi_1^{\delta}(y-y_1)\varphi_2^{\delta}(\xi-\xi_1)\:dy_1\:d\xi_1\\
&&+ \int\left[a_{N+1}(y,\xi)-a_{N+1}(y_1,\xi_1)\right]\psi(y_1,\xi_1)\varphi_1^{\delta}(y-y_1)\p_{xi}\varphi_2^{\delta}(\xi-\xi_1)\:dy_1\:d\xi_1\
\end{eqnarray*}

Thus, we compute, for $(y,y_1,\xi,\xi_1)\in \R^{2N+2}$, $1\leq
i\leq N+1$,
\begin{eqnarray*}
a_i(y,\xi)-a_i(y_1,\xi_1)&=&(y-y_1)\cdot\int_0^1\nabla_y a_i(\tau
y + (1-\tau) y_1, \tau \xi + (1-\tau) \xi_1)\:d\tau\\
&&+(\xi-\xi_1)\cdot\int_0^1\p_{\xi} a_i(\tau y + (1-\tau) y_1,
\tau \xi + (1-\tau) \xi_1)\:d\tau.
\end{eqnarray*}
Set, for $1\leq k,i\leq N$, $y\in\R^N$, $\xi\in\R$,
\begin{gather*}
 \phi_{k,i}(y,\xi)=y_k \frac{\p \varphi_1}{\p y_i}(y)\varphi_2(\xi),\qquad
\phi_{k,N+1}(y,\xi)=y_k\frac{\p \varphi_2}{\p \xi}(\xi)\varphi_1(y),\\
\zeta_i(y,\xi)=\xi\frac{\p \varphi_1}{\p y_i}(y)\varphi_2(\xi),\qquad
\zeta_{N+1}(y,\xi)=\xi\frac{\p \varphi_2}{\p \xi}(\xi)\varphi_1(y).
\end{gather*}

Notice that
$$
 \int_{\R^{N+1}}\phi_{k,i}=-\delta_{k,i},\qquad
\int_{\R^{N+1}}\zeta_i =-\delta_{N+1,i}.$$

Then
\begin{eqnarray*}
r_{\delta}(y,\xi)&=&\sum_{i=1}^{N+1}\sum_{k=1}^N\int\frac{\p a_i}{\p y_k}(\tau y + (1-\tau) y_1, \tau \xi + (1-\tau)
\xi_1)\psi(y_1,\xi_1)\phi_{k,i}^{\delta}(y-y_1,\xi-\xi_1)\:dy_1\:d\xi_1\:d\tau\\
&& + \sum_{i=1}^{N+1}\int\frac{\p a_i}{\p \xi}(\tau y + (1-\tau) y_1, \tau \xi + (1-\tau)
\xi_1)\psi(y_1,\xi_1)\zeta_i^{\delta}(y-y_1,\xi-\xi_1)\:dy_1\:d\xi_1\:d\tau.
\end{eqnarray*}
Hence as $\delta \to 0$, $r_{\delta}$ converges to
$$-\dv_{y,\xi}(a(y,\xi))\:\psi(y,\xi)=0$$
in $L^p_{\text{loc}}(\R^{N+1})$ for any $p<\infty$ and for all
$(t,x)\in [0,\infty)\times\R^N$. We now pass to the limit as $\delta\to 0$, with $\e$ fixed, and we obtain
\begin{eqnarray*}
&&-\int \fe(s,z,\xi)\left\{\p_t \varphi(t-s,x-z) +  a_i\left(\frac{z}{\e},\xi\right) \p_{i} \varphi(t-s,x-z)\right\}\psi\left(\frac{z}{\e},\xi\right) \theta (t,x)\:d\xi\:dx\:dz\:ds\:dt\\
&&- \alpha_- \int\theta(t,x)\:\p_{\xi} K(\xi)A\left( \frac{z}{\e},\xi \right)\cdot \nabla_x \varphi(t-s,x-z) \: dt\:dx \:ds\:dz\:d\xi \\
&\geq &0.
\end{eqnarray*}
Passing to the limit as $\e$ vanishes, we are led to
\begin{eqnarray*}
&&-\int \!f^0(s,z,y,\xi)\left\{\p_t \varphi(t-s,x-z) +  a_i\left(y,\xi\right) \p_{i} \varphi(t-s,x-z)\right\}\psi(y,\xi) \theta (t,x)\:d\xi\:dx\:dz\:ds\:dy\:dt\\
&&- \alpha_- \int\theta(t,x)\:\p_{\xi} K(\xi)A\left( y,\xi \right)\cdot \nabla_x \varphi(t-s,x-z) \: dt\:dx \:ds\:dy\:dz\:d\xi \\
&\geq &0.
\end{eqnarray*}
Since
$$
\int\theta(t,x)\nabla_x \varphi(t-s,x-z) \: dt\:dx \:ds\:dz=-\left(\int\theta(t,x)\:dt \:dx\right) \left( \int \nabla_z \varphi(s,z)\:ds\:dz \right)=0,$$
we deduce that
$$\int\! f^0(s,z,y,\xi)\left\{\p_t \varphi(t-s,x-z) + \Si a_i\left(y,\xi\right) \p_{i} \varphi(t-s,x-z)\right\}\psi\left(y,\xi\right) \theta (t,x)d\xi\:dx\:dz\:ds\:dy\:dt\leq 0,$$
which means that $f^0$ satisfies condition \eqref{hyp:M}. There only remains to check the strong continuity of $f$ at time $t=0$.

\subsection{Strong continuity at time $t=0$}

The continuity property for$ f^0$ is inherited from uniform continuity properties at
time $t=0$ for the sequence $\fe$. This is strongly linked to the well-preparedness of the initial data (condition \eqref{eq:cell_kin}), that is, the fact that for all $x\in\R^N$, $u_0(x,\cdot)$ is an entropy solution of the cell problem
$$
\dv_y A(y,u_0(x,y))=0.
$$

First, let us consider a regularization of the initial data
$$
g_n^{\delta}=f_0\ast_x \rho_n \ast_y\varphi_1^{\delta}\ast_{\xi}\varphi_2^{\delta}.
$$
with $\rho_n $ a convolution kernel ($n\in\N$), $\delta>0$, and $\varphi_i^{\delta}$ defined as in the previous subsection.
Then we can write
\begin{eqnarray}
 &&\Si a_i\left(\frac x \e,\xi \right)\cdot\frac{\p}{\p x_i} \left[g_n^{\delta}\left(x,\frac x \e,\xi
\right) \right]+ \frac{1}{\e} a_{N+1}\left(\frac x \e,\xi \right)\frac{\p}{\p \xi} g_n^{\delta}\left(x,\frac x \e,\xi
\right)\nonumber \\
&=& \frac{1}{\e}a\left(\frac x \e,\xi
\right)\cdot \left(\nabla_{y,\xi} g_n^{\delta}\right)\left(x,\frac{x}{\e} ,\xi \right) + \Si a_i\left(\frac x \e,\xi \right)\left(\frac{\p}{\p x_i} g_n^{\delta} \right)\left(x,\frac x \e,\xi
\right)\label{eq:gn}\\
&:=&r^{\e}_{n,\delta}.\nonumber
\end{eqnarray}

Notice that
$$
||\nabla_x g_n^{\delta}||_{L^\infty(\R^N\times Y\times \R)}\leq ||\nabla_x \rho_n||_{L¹(\R^N)},
$$
and $$a\left(y,\xi
\right)\nabla_y g_n^{\delta}\left(x,y ,\xi \right)=\p_{\xi } m_n^\delta + r_n^\delta, $$
where
\begin{gather*}
 m_n^\delta= m_0 \ast_x \rho_n \ast_y\varphi_1^{\delta}\ast_{\xi}\varphi_2^{\delta},\\
 r_n^\delta(x,y,\xi) =a\left(y,\xi
\right)\nabla_y g_n^{\delta}\left(x,y ,\xi \right) - \left[a f_0\ast_x \rho_n  \right] \ast_{y,\xi} \nabla_{y,\xi} \varphi_1^\delta(y)\varphi_2^\delta(\xi).
\end{gather*}
Then for all $n\in \N$, for all $x\in\R^N$, $r_n^\delta$ vanishes as $\delta \to 0$ in $L^1_{\text{loc}}(Y\times \R)$ and almost everywhere. The proof of this fact is exactly the same as in the preceding subsection, and thus, we leave the details to the reader. As a consequence,
$$
r^{\e}_{n,\delta}(x,\xi)= \frac{1}{\e}\p_{\xi} m_n^\delta\left(x,\frac{x}{\e} ,\xi \right) + R^{\e}_{n,\delta}(x,\xi),
$$
and there exists a constant $C_n$, independent of $\e$, such that for all $n\in\N$, for all $\e>0$, and for almost every $x,\xi$
$$
\limsup_{\delta \to 0}|R^{\e}_{n,\delta}(x,\xi)|\leq C_n.
$$
Moreover, $\supp\: R^{\e}_{n,\delta}\subset \R^N\times [-R-1,R+1]$ if $\delta<1$.

Now, we multiply \eqref{eq:depart_kf_gal} by $1 - 2g_n^\delta\left(x,x/\e,\xi \right)$, and \eqref{eq:gn} by $1- 2 \fe(t,x,\xi) $. Setting
\begin{eqnarray*}
 h_{n,\delta}^\e(t,x,\xi)&:= &\fe(t,x,\xi) \left[ 1 - 2g_n^\delta\left(x,\frac{x}{\e},\xi \right) \right] + g_n^\delta\left(x,\frac{x}{\e},\xi \right)\left[1- 2 \fe(t,x,\xi)\right]\\
&=&\left|\fe(t,x,\xi) -  g_n^\delta\left(x,\frac{x}{\e},\xi \right) \right|^2 +  g_n^\delta\left(x,\frac{x}{\e},\xi \right) - \left|g_n^\delta\left(x,\frac{x}{\e},\xi \right) \right|^2,
\end{eqnarray*}
we obtain
\begin{multline}
\frac{\p}{\p t} h_{n,\delta}^\e(t,x,\xi)+\Si a_i\left(\frac x \e,\xi
\right)\p_{x_i} h_{n,\delta}^\e(t,x,\xi) + \frac{1}{\e} a_{N+1}\left(\frac x \e,\xi
\right)\p_{\xi} h_{n,\delta}^\e(t,x,\xi)=\\=\frac{\p \me }{\p \xi}\left[1 -
2g_n^{\delta}\left(x,\frac x \e,\xi \right)\right] +
\frac{1}{\e}\p_{\xi} m_n^\delta\left(x,\frac{x}{\e} ,\xi \right) \left[1-2
\fe(t,x,\xi) \right] + R^{\e}_{n,\delta}(x,\xi) \left[1-2
\fe(t,x,\xi) \right] .\label{eq:comp_fe_f0} \end{multline}
Notice that
\begin{gather*}
 \p_{\xi}\left[1-2
\fe(t,x,\xi) \right]=2 \delta(\xi=\ue(t,x)),\\
\frac{\p}{\p\xi}\left(1 - 2g_n^\delta\left(x,\frac x \e,\xi
\right)\right)=2 \nu_{n,\e,\delta}(x,\xi),
\end{gather*}
where $\nu_{n,\e,\delta}$ is a non-negative function in $\mathcal  C^\infty(\R^{N+1})$, with support in $\R^N\times [-M-1,M+1]$ if $\delta<1$. Notice also that $\fe(t,x,\xi)-g_n^{\delta}\left(x,x/\e,\xi
\right)=0$ if $|\xi|$ is large enough ($|\xi|>M+1$). Take a cut-off function $\zeta=\zeta(x)$ such that
$\zeta(x)=e^{-|x|}$ when $|x|\geq 1$, and $\frac 1 e \leq
\zeta(x)\leq 1$ for $|x|\leq 1$. Then there exists a
constant $C$ such that
$$
|\nabla_x \zeta(x)|\leq C \zeta(x)\quad
\forall x\in\R^N.
$$

Hence, mutliplying \eqref{eq:comp_fe_f0} by $\zeta(x)$ and integrating on $\R^{N+1}$, we obtain a bound of the type
\begin{eqnarray*}
\frac{d}{dt} \int_{\R^{N+1}} h_{n,\delta}^\e(t,x,\xi)\zeta(x)\:dx\:d\xi &\leq & C \int_{\R^{N+1}} h_{n,\delta}^\e(t,x,\xi)\zeta(x)\:dx\:d\xi\\&&+   \int_{\R^{N+1}}\left|R^{\e}_{n,\delta}(x,\xi) \right|\;\left| 1-2
\fe(t,x,\xi) \right|\zeta(x)\:dx\:d\xi.
\end{eqnarray*}
Using Gronwall's lemma and passing to the limit as $\delta\to 0$ with $\e$ and $n\in\N$ fixed, we retrieve, for all $t\geq 0$,
\begin{eqnarray*}
\int_{\R^{N+1}}\left|\fe(t,x,\xi)-g_n\left(x,\frac x
\e,\xi \right)\right|^2\zeta(x)\:dx\:d\xi&\leq & e^{Ct}\int_{\R^{N+1}}\left|f_0\left(x,\frac{x}{\e} ,\xi \right)-g_n\left(x,\frac x
\e,\xi \right)\right|^2\zeta(x)\:dx\:d\xi\\
&+&  e^{Ct}\int_{\R^{N+1}}\!\!\left[g_n\left(x,\frac x
\e,\xi \right)-\left|g_n\left(x,\frac x
\e,\xi \right)\right|^2\right]\zeta(x)\:dx\:d\xi\\
&+&   C_n(e^{Ct}- 1),
\end{eqnarray*}
where the constant $C_n$ does not depend on $\e$, and $g_n= f_0\ast_x \rho_n$. And for all $n\in\N$, $\e>0$, we have
\begin{eqnarray*}
&&\int_{\R^{N+1}}\left|f_0\left(x,\frac{x}{\e} ,\xi \right)-g_n\left(x,\frac x
\e,\xi \right)\right|^2\zeta(x)\:dx\:d\xi\\
&\leq & \int_{\R^{N+1}}\int_{\R^N}\left|f_0\left(x,\frac{x}{\e} ,\xi \right) - f_0\left(x',\frac{x}{\e} ,\xi \right) \right|^2\rho_n(x-x')\zeta(x)\:dx\:dx'\:d\xi\\
&\leq & \int_{\R^{N}}\int_{\R^N}\left|u_0\left(x,\frac{x}{\e} ,\xi \right) - u_0\left(x',\frac{x}{\e} ,\xi \right) \right|\rho_n(x-x')\zeta(x)\:dx\:dx'\\
&\leq &  \int_{\R^{N}}\int_{\R^N}\sup_{y\in Y}\left|u_0\left(x,y ,\xi \right) - u_0\left(x',y ,\xi \right) \right|\rho_n(x-x')\zeta(x)\:dx\:dx'.
\end{eqnarray*}
The right-hand side of the above inequality vanishes as $n\to\infty$ because $u_0\in L^1_{\text{loc}}(\R^N, \mathcal C_{\text{per}}(Y))$. Similarly,
\begin{eqnarray*}
&&\int_{\R^{N+1}}\left[g_n\left(x,\frac x
\e,\xi \right)-\left|g_n\left(x,\frac x
\e,\xi \right)\right|^2\right]\zeta(x)\:dx\:d\xi\\
&\leq & \int_{\R^{N+1}}\left[g_n\left(x,\frac x
\e,\xi \right)-f_0\left(x,\frac x
\e,\xi \right)\right]\zeta(x)\:dx\:d\xi\\
& &+ \int_{\R^{N+1}}\left[f_0\left(x,\frac x
\e,\xi \right)^2-g_n\left(x,\frac x
\e,\xi \right)^2\right]\zeta(x)\:dx\:d\xi\\
&\leq & 3 \int_{\R^{N+1}}\left|g_n\left(x,\frac x
\e,\xi \right)-f_0\left(x,\frac x
\e,\xi \right)\right|\zeta(x)\:dx\:d\xi\\
&\leq & 3 \int_{\R^{N}}\int_{\R^N}\sup_{y\in Y}\left|u_0\left(x,y ,\xi \right) - u_0\left(x',y ,\xi \right) \right|\rho_n(x-x')\zeta(x)\:dx\:dx'.
\end{eqnarray*}

Hence, we deduce that there exists a function
$\omega:[0,\infty)\to [0,\infty)$, independent of $\e$ and satisfying $\lim_{t\to 0} \omega(t)=0$, such that
$$
\int_{\R^{N+1}}\left|\fe(t,x,\xi)-f_0\left(x,\frac{x}{\e}  ,\xi\right)\right|\zeta(x)\:dx\:d\xi\leq \omega(t)
$$
for all $t>0$.

Then, we prove that the same property holds for the function
$f^0$. Indeed, we write
$$
\left|\fe(t,x,\xi)-\mathbf 1_{\xi< u_0\left( x,\frac{x}{\e} \right)}\right|=\fe - 2
\fe\mathbf 1_{\xi< u_0\left( x,\frac{x}{\e} \right)} +
\mathbf 1_{\xi< u_0\left( x,\frac{x}{\e} \right)};
$$
let $\theta\in L^\infty ([0,\infty))$ with compact support and such that $\theta\geq 0$. Then for all $\e>0$,
$$
\int_0^\infty\int_{\R^{N+1}}\left[\fe - 2
\fe\mathbf 1_{\xi< u_0\left( x,\frac{x}{\e} \right)} +
\mathbf 1_{\xi< u_0\left( x,\frac{x}{\e} \right)}\right]\zeta(x)\theta(t)\:dx\:d\xi\:dt\leq \int_0^\infty\omega(t)\theta(t)\:dt.
$$
Since $u_0\in L^1_{\text{loc}}(\R^N, \mathcal C_{\text{per}}(Y))$, it is an
admissible test function in the sense of G. Allaire (see
\cite{Allaire}); we deduce that $\mathbf 1_{\xi<u_0}$ is also an admissible test function. This is not entirely obvious because it is a discontinuous function of $u_0$. However, this difficulty can be overcome thanks to an argument similar to the one developed below in subsection \ref{ssec:stcv}, and which we do not reproduce here. Thus, we can pass to the two-scale limit in the
above inequality. We obtain
\begin{multline*}
 \int_0^{\infty}\int_{\R^{N+1}\times Y}(f^0(t,x,y,\xi)-|f^0(t,x,y,\xi)|^2 +
|f^0(t,x,y,\xi)
-\mathbf 1_{\xi< u_0\left( x,y\right)}|^2\theta(t)\zeta(x)\:dt\:dx\:dy\:d\xi\leq\\
\leq
\int_0^{\infty}\theta(t)
\omega(t)\:dt
\end{multline*}
Notice that $|f^0|-|f^0|^2 \geq 0$ almost everywhere. As a consequence, taking $\theta (t)=\mathbf 1_{0<t<\tau}$, with $\tau>0$ arbitrary, we deduce that
$$
\frac{1}{\tau}\int_0^\tau |f^0(t)
-\chi(\xi,u_0(x,y))|^2\zeta(x)\:dt\:dx\:dy \leq \frac 1 \tau \int_0^\tau \omega(t)\:dt,
$$
and the left-hand side vanishes as $\tau\to 0$. Thus the continuity property is satisfied at time $t=0$.

\vskip3mm

Hence, we have proved that any two-scale limit of the sequence $\fe$ is a solution of the limit system. Thus the existence result in Theorem \ref{thm:lim_ev_eq} is proved, as well as the convergence result of Theorem \ref{thm:strg_cvg}. We now tackle the proof of the uniqueness and rigidity results of Theorem \ref{thm:lim_ev_eq}. The strong convergence result of Theorem \ref{thm:homogenized_sep} will follow from the rigidity.

\section{Uniqueness of solutions of the limit evolution problem\label{sec:contraction}}

In this section, we prove the second and the third point in Theorem \ref{thm:lim_ev_eq}, that is, if $f$ is any solution of the limit evolution problem, then there exists a function $u\in L^{\infty}([0,\infty)\times \R^N \times Y)$ such that $f(t,x,y,\xi)=\mathbf 1_{\xi<u(t,x,y)}$ almost everywhere, and if $f_1=\mathbf 1_{\xi<u_1}$, $f_2=\mathbf 1_{\xi<u_2}$ are two generalized kinetic solutions, then the contraction principle \eqref{in:contraction} holds.

\subsection{The rigidity result \label{ssec:rigidity}}

Let $f$ be a generalized kinetic solution of the limit problem, with initial data $\mathbf 1_{\xi<u_0}$. The rigidity result relies on the comparison between $f$ and $f^2$. Precisely, we prove that $f=|f|^2$ almost everywhere, and since $\p_{\xi} f=  -\nu\leq 0$, this identity entails that there exists a function $u$ such that $f=\mathbf 1_{\xi<u}$. Thus, we now turn to the derivation of the equality $|f|=|f|^2$.

Let $\delta>0$ arbitrary, and let $\theta_1\in\mathcal D(\R), \theta_2\in\mathcal D(\R^N)$ such that
\begin{gather*}
\theta_1\geq 0, \ \theta_2\geq 0 ,\\
\int_{\R}\theta_1=\int_{\R^N}\theta_2=1,\\
\supp\: \theta_1\subset [-1,0]\text{ and } \theta_1(0)=0.
\end{gather*}
We set, for $(t,x)\in\R^{N+1}$
$$
\theta^{\delta}(t,x)=\frac{1}{\delta^{N+1}}\theta_1\left ( \frac{t}{\delta}\right)\theta_2\left ( \frac{x}{\delta}\right).
$$
Set $f^\delta:=f\ast_{t,x}\theta^{\delta}$, $\mathcal M^\delta:=\mathcal M\ast_{t,x}\theta^{\delta}$. Then $f^{\delta}$ is a solution of
$$
\frac{\p f^\delta}{\p t}+ \Si a_i(y,\xi)\frac{\p f^\delta}{\p x_i} = \mathcal M^\delta.
$$
Moreover, $f^\delta$ satisfies the following properties
\begin{gather}
0\leq f^\delta\leq 1,\label{pro:fdelta1}\\
\dv_{y,\xi}(a(y,\xi) f^\delta)=\p_{\xi} m\ast_{t,x}\theta^{\delta},\\
\p_{\xi }f^\delta=-\nu\ast_{t,x}\theta^{\delta},
\end{gather}
whereas $\mathcal M^\delta$ satisfies
\begin{gather}
\mathcal M^\delta\in\mathcal C((0,T)\times \R^N, L^2(Y\times \R))\cap L^\infty([0,\infty)\times \R^N \times Y\times \R),\\
\mathcal M^\delta(\cdot,\xi)=0\quad \text{if }|\xi|>M,
f^\delta(\cdot,\xi)=0\quad \text{if }\xi>M,\quad f^\delta(\cdot,\xi)=1\quad \text{if }\xi<-M,
\\
\int_{ Y\times \R}\mathcal M^\delta \psi \leq 0\quad \forall \psi \in\mathcal G.\label{pro:fdelta2}
\end{gather}
In particular, notice that $(1-2f^\delta(t,x))\in\mathcal G$ for all $t,x$, and $f^\delta(t,x,y,\xi) - f^\delta(t,x,y,\xi)^2=0$ if $|\xi|>M$.

Let $\zeta\in\mathcal C^{\infty}(\R^N)$ be a cut-off function as in the previous subsection. We multiply the equation on $f^\delta $ by $(1 - 2 f^\delta ) \zeta(x)$, and we integrate over $\R^N\times Y \times \R$. We obtain
$$
\frac{d}{dt} \int_{\R^N\times Y \times\R} \left(f^\delta -|f^\delta |^2 \right)\zeta-  \int_{\R^N\times Y \times\R}a_i(y,\xi) \p_i \zeta(x) \left(f^\delta -|f^\delta |^2 \right)= \int_{\R^N\times Y \times\R}\mathcal M^{\delta} \left(1- 2 f^\delta  \right)\zeta\leq 0.
$$
We then deduce successively, using Gronwall's lemma,
\begin{gather}
 \frac{d}{dt} \int_{\R^N\times Y \times\R} \left(f^\delta -|f^\delta |^2 \right)\zeta\leq C\int_{\R^N\times Y \times\R} \left(f^\delta -|f^\delta |^2 \right)\zeta,\nonumber\\
 \int_{\R^N\times Y \times\R} \left(f^\delta(t) -|f^\delta(t) |^2 \right)\zeta\leq e^{Ct}\int_{\R^N\times Y \times\R} \left(f^\delta (t=0)-|f^\delta (t=0)|^2 \right)\zeta\quad \forall t>0,\nonumber\\
\int_0^T \int_{\R^N\times Y \times\R} \left(f^\delta-|f^\delta |^2 \right)\zeta\leq\frac{e^{CT}-1}{C}\int_{\R^N\times Y \times\R} \left(f^\delta (t=0)-|f^\delta (t=0)|^2\right)\zeta,\label{est:f_f2_rigid}
\end{gather}
with a constant $C$ depending only on $||a||_{L^\infty(Y\times [-R,R])}$.

Let us now check that $f^\delta(t=0)$ strongly converges towards $\mathbf 1_{\xi<u_0}=f_0$ at time $t=0$. In fact, the main difference between the proof of Theorem \ref{thm:lim_ev_eq} and the one for generalized kinetic solutions of scalar conservation laws (see chapter 4 in \cite{BP}) lies in this particular point. Indeed, in the case of scalar conservation laws, the continuity property can be inferred from the equation itself; in the present case, the lack of structure of the right-hand side $\mathcal M$ prevents us from deriving such a result, and hence the continuity  of solutions at time $t=0$ is a necessary assumption in definition \ref{def:lim_ev_eq}.

Using hypothesis \eqref{hyp:cont_time0}, we write, for almost every $x,y,\xi$,
\begin{eqnarray*}
 f^\delta(t=0,x,y,\xi)&=&\int_{\R^{N+1}}f(s,z,y,\xi)\theta^{\delta}(-s,x-z)\:ds\:dz\\
 f^\delta(t=0,x,y,\xi)- f_0\ast_x  \theta_2^{\delta}(x,y,\xi)&=&\int_{\R^{N+1}}\left(f(s,z,y,\xi)- f_0(z,y,\xi)\right)\theta^{\delta}(-s,x-z)\:ds\:dz.\end{eqnarray*}
As a consequence, for all $\delta>0$
\begin{eqnarray*}
&&\int_{\R^N\times Y\times \R}\left|f^\delta(t=0) - f_0\ast_x  \theta_2^{\delta}\right|^2\zeta(x)\:dx\:dy\:d\xi \\
&\leq & \int_{\R^N\times Y\times \R}\int_{\R^{N+1}}\left|f(s,z,y,\xi)- f_0(z,y,\xi)\right|^2\zeta(x)\theta^{\delta}(-s,x-z)\:dx\:dy\:d\xi\:ds\:dz\\
&\leq & \int_\R||f(s)-f_0||_{L^2(\R^N\times Y\times \R, \zeta(x)\:dx\:dy\:d\xi)}^2\frac{1}{\delta}\theta_1\left ( \frac{-s}{\delta}\right)\:ds\:dx\:dy\:d\xi
+ 2R |Y|\:||\zeta - \zeta\ast\check{\theta}_2^{\delta} ||_{L^1(\R^N)}\\
&\leq & \frac{C}{\delta}\int_0^{\delta}||f(s)-f_0||_{L^2(\R^N\times Y\times \R, \zeta(x)\:dx\:dy\:d\xi)}^2\:ds+ 2R |Y|\:||\zeta - \zeta\ast\check{\theta}_2^{\delta} ||_{L^1(\R^N)}.
\end{eqnarray*}
The right-hand side of the last inequality vanishes as $\delta\to 0$, and thus $ f^\delta(t=0)$ converges towards $f_0$ as $\delta\to 0$ in $L^2(\R^N\times Y\times \R, \zeta(x)\:dx\:dy\:d\xi)$, and hence also in $L^1(\R^N\times Y\times \R, \zeta(x)\:dx\:dy\:d\xi)$.
Consequently,
$$
\int_{\R^N\times Y \times\R} \left(f^\delta (t=0)-|f^\delta (t=0)|^2\right)\zeta \to 0\quad \text{as }\delta\to 0.
$$
Above, we have used the fact that $f_0=\mathbf 1_{\xi<u_0}$, and thus $f_0= f_0^2$.

Now, we pass to the limit as $\delta\to 0$ in \eqref{est:f_f2_rigid}; we obtain, for all $T>0$,
$$
\int_0^T \int_{\R^N\times Y \times\R} \left(f-|f|^2 \right)\varphi \leq 0.
$$
Since the integrand in the left-hand side is non-negative, we deduce that $|f|=|f|^2$ almost everywhere. The rigidity property follows.

\subsection{Contraction principle\label{ssec:contrac}}

Let $f_1$, $f_2$ be two generalized kinetic solutions of the limit problem; we denote by $M_1,M_2$, and $\mathcal M_1, \mathcal M_2$, the constants and distributions associated to $f_1$, $f_2$, respectively. Without loss of generality, we assume that $M_1\leq M_2$. According to the rigidity result, there exist functions $u_1,u_2\in L^\infty([0,\infty)\times\R^N\times Y)\cap L^\infty([0,\infty), L^1(\R^N\times Y))$ such that $f_i=\mathbf 1_{\xi<u_i}$.

As in the previous subsection, we regularize $f_i, \mathcal M_i$ by convolution in the variables $t,x$, and we denote by $f_i^{\delta}, \mathcal M_i^{\delta}$ the functions thus obtained. The strategy of the proof is the same as in \cite{BP}, Theorem 4.3.1. The idea is to derive an inequality of the type
\be
\frac {d} {dt} \int|f_1(t,x,y,\xi)-f_2(t,x,y,\xi)|\zeta(x)\:dx\:dy\:d\xi\leq C\int|f_1(t,x,y,\xi)-f_2(t,x,y,\xi)|\zeta(x)\:dx\:dy\:d\xi ,
\ee
where $\zeta$ is a cut-off function as in the previous section.

Since $|f_1(t)-f_2(t)|=|f_1(t)-f_2(t)|²=f_1+ f_2- 2 f_1f_2$, let us first write the equation satisfied by $g^{\delta}:=f_1^{\delta}+ f_2^{\delta}- 2 f_1^{\delta}f_2^{\delta}$. We compute
\begin{gather*}
\left\{\p_ t  f_1^{\delta} + \Si a_i(y,\xi)\frac{\p}{\p x_i}f_1^{\delta}= \mathcal M_1^{\delta}\right\}\quad \times 1 - 2 f_2^{\delta},\\
\left\{\p_ t  f_2^{\delta} + \Si a_i(y,\xi)\frac{\p}{\p x_i}f_2^{\delta}= \mathcal M_2^{\delta}\right\}\quad \times 1 - 2 f_1^{\delta}.
\end{gather*}
Adding the two equations thus obtained leads to
$$
\p_ t  g^{\delta} + \Si a_i(y,\xi)\frac{\p}{\p x_i} g^{\delta}=\mathcal M_1^{\delta}\left[1 - 2 f_2^{\delta}\right] + \mathcal M_2^{\delta}\left[1 - 2 f_1^{\delta}\right].
$$
Notice that thanks to \eqref{hyp:compact_support2}, \eqref{hyp:compact_support3} and the microscopic constraints \eqref{eq:cell_kin}, \eqref{hyp:dfdxi}, $1 - 2 f_i^{\delta}(t,x)\in\mathcal G$ for all $(t,x)$. Hence
$$
\int_{Y\times \R}\mathcal M_2^{\delta}(t,x)\left[1 - 2 f_1^{\delta}(t,x)\right]\leq 0\quad \forall (t,x)\in[0,\infty)\times \R^N,
$$
and the same inequality holds if the roles of $f_1$ and $f_2$ are exchanged.

Now, take a cut-off function $ \zeta \in\mathcal C^{\infty}(\R^N)$ satisfying the same assumptions as in the previous subsection; multiply the equation on $g^\delta$ by $\zeta( x)$, and integrate over $\R^N\times Y\times \R$; this yields
$$
\frac d {dt} \int_{\R^N\times Y\times \R}g^\delta(t,x,y,\xi) \zeta( x)\:dx\:dy\:d\xi \leq C\int_{\R^N\times Y\times \R}g^\delta(t,x,y,\xi) \zeta( x)\:dx\:dy\:d\xi \quad \forall t>0,
$$
and thus
$$
 \int_{\R^N\times Y\times \R}g^\delta(t,x,y,\xi) \zeta( x)\:dx\:dy\:d\xi \leq e^{Ct} \int_{\R^N\times Y\times \R}g^\delta(t=0,x,y,\xi) \zeta(x)\:dx\:dy\:d\xi.
$$
According to the strong convergence results of $f_i^{\delta}(t=0)$ derived in the previous section, we can pass to the limit as $\delta\to 0$. We infer that for almost every $t>0$,
\begin{eqnarray}
&& \int_{\R^N\times Y\times \R}|f_1(t,x,y,\xi)-f_2(t,x,y,\xi)| \zeta( x)\:dx\:dy\:d\xi \nonumber\\
&\leq & e^{Ct} \int_{\R^N\times Y\times \R}|f_1(t=0,x,y,\xi)-f_2(t=0,x,y,\xi)| \zeta( x)\:dx\:dy\:d\xi.\label{in:contrac_zeta}
\end{eqnarray}

This completes the derivation of the contraction principle for the limit system. Uniqueness of solutions of the limit system follows. In particular, we deduce that the whole sequence $\fe$ of solutions of \eqref{eq:depart_kf_gal} two-scale converges towards $f^0$.

\subsection{Strong convergence result\label{ssec:stcv}}

Here, we explain why the strong convergence result stated in Theorem \ref{thm:strg_cvg} holds, that is, we prove \eqref{strongCVregul}. This fact is rather classical, and is a direct consequence of the fact that
$$
\mathbf 1_{\xi<\ue(t,x)}\stackrel{\text{2 sc.}}{\rightharpoonup} \mathbf 1_{\xi<u(t,x,y)}.
$$
Let us express this result in terms of Young measures: the above two-scale convergence is strictly equivalent to the fact that the two-scale Young measure $\nu_{t,x,y}$ associated with the sequence $\ue$ is the Dirac mass $\delta (\xi=u(t,x,y))$ (see \cite{BP}, Chapter 2). And it is well-known (see \cite{wetransport}) that if $u$ is a smooth function, then
$$
d\nu_{t,x,y}(\xi)=\delta (\xi=u(t,x,y))\quad \iff \quad \ue-u\left( t,x,\frac{x}{\e} \right)\to 0\quad \text{in }L^1_{\text{loc}}.$$

For the reader's convenience, we now prove the result without using two-scale Young measures. We define $u_{\delta}= u\ast_x \varphi_{\delta}$, with $ \varphi_{\delta}$ a standard mollifier. Take $K\in\mathcal D(\R)$ such that $0\leq K\leq 1$, and $K(\xi)=1$ if $|\xi|\leq M$. Take also a sequence $\theta_n\in\mathcal C^\infty(\R)$ such that $0\leq \theta_n\leq 1$, and
$$
\theta_n(\xi)= 1 \text{ if } \xi <-\frac{1}{n},\quad \theta_n(\xi)= 0\text{ if } \xi >\frac{1}{n}.
$$

Then we write
\begin{eqnarray*}
\left|\mathbf 1_{\xi<\ue(t,x)} - \mathbf 1_{\xi< u_{\delta}\left(t,x,\frac{x}{\e}  \right)}\right|^2&=& \mathbf 1_{\xi<\ue(t,x)} - 2 \mathbf 1_{\xi< u_{\delta}\left(t,x,\frac{x}{\e}\right)}\mathbf 1_{\xi<\ue(t,x)} + \mathbf 1_{\xi< u_{\delta}\left(t,x,\frac{x}{\e}\right)} \\
&=&\mathbf 1_{\min \left( \ue(t,x) , u_{\delta}\left(t,x,\frac{x}{\e}  \right) \right)<\xi<\max  \left( \ue(t,x) , u_{\delta}\left(t,x,\frac{x}{\e}  \right) \right)}.
\end{eqnarray*}
The function $ \mathbf 1_{\xi< u_{\delta}\left(t,x,\frac{x}{\e}  \right)}$ is not smooth enough to be used as an oscillating test function. Thus we replace it by
$$
\theta_n\left(\xi -u_{\delta}\left(t,x,\frac{x}{\e}  \right)\right),
$$
and we evaluate the difference : for all compact set $C\subset [0,\infty)\times \R^N$,
$$
\int_C \int_{\R}\left| \mathbf 1_{\xi< u_{\delta}\left(t,x,\frac{x}{\e}  \right)} -  \theta_n\left(\xi -u_{\delta}\left(t,x,\frac{x}{\e}  \right)\right)\right|K(\xi)\:dt\:dx \: d\xi \leq \frac{2}{n } |C|.
$$
According to the two-scale convergence result, for all $n\in\N$,
\begin{multline*}
 \int_C \int_{\R}\theta_n\left(\xi -u_{\delta}\left(t,x,\frac{x}{\e}  \right)\right)\mathbf 1_{\xi<\ue(t,x)}K(\xi)\:dt\:dx\:d\xi\to \\
\to\int_C \int_{\R}\theta_n\left(\xi -u_{\delta}\left(t,x,y  \right)\right)\mathbf 1_{\xi<u(t,x,y)}K(\xi)\:dt\:dx\:dy\:d\xi.
\end{multline*}
Since the sequence $\theta_n\left(\xi -u_{\delta}\right) $ uniformly converges towards $\mathbf 1_{\xi<u_\delta}$ as $n\to\infty$, we can pass to the limit as $n\to\infty$, and we deduce
$$
\int_C \int_{\R} \mathbf 1_{\xi< u_{\delta}\left(t,x,\frac{x}{\e}  \right)}\mathbf 1_{\xi<\ue(t,x)}K(\xi)\:dt\:dx\:d\xi\to \int_C \int_{\R\times Y} \mathbf 1_{\xi< u_{\delta}\left(t,x,y  \right)}\mathbf 1_{\xi<u(t,x,y)}K(\xi)\:dt\:dx\:dy\:d\xi.
$$
Simlarly, as $\e\to 0$, for all $\delta>0$,
\begin{gather*}
 \int_C \int_{\R} \mathbf 1_{\xi< u_{\delta}\left(t,x,\frac{x}{\e}  \right)}K(\xi)\:dt\:dx\:d\xi\to \int_C \int_{\R\times Y}\mathbf 1_{\xi< u_{\delta}\left(t,x,y  \right)}K(\xi)\:dt\:dx\:dy\:d\xi,\\
\int_C \int_{\R\times Y} \mathbf 1_{\xi< \ue(t,x)}K(\xi)\:dt\:dx\:d\xi\to \int_C \int_{\R}\mathbf 1_{\xi<u\left(t,x,y  \right)}K(\xi)\:dt\:dx\:dy\:d\xi.
\end{gather*}
Thus
$$
\int_C \int_{\R} \left|\mathbf 1_{\xi<\ue(t,x)} - \mathbf 1_{\xi< u_{\delta}\left(t,x,\frac{x}{\e}  \right)}\right|^2 K(\xi)\:dt\:dx\:d\xi \to \int_C \int_{\R\times Y} \left|\mathbf 1_{\xi<u\left(t,x,y  \right)} - \mathbf 1_{\xi<u_{\delta}\left(t,x,y\right)}\right| \:K(\xi)\:dt\:dx\:dy\:d\xi
$$

On the other hand,
\begin{gather*}
 \int_C \int_{\R} \left|\mathbf 1_{\xi<\ue(t,x)} - \mathbf 1_{\xi< u_{\delta}\left(t,x,\frac{x}{\e}  \right)}\right|^2 K(\xi)\:dt\:dx\:d\xi= \left|\left| \ue(t,x) - u_{\delta} \left(t,x,\frac{x}{\e}  \right)\right|  \right|_{L^1(C)},\\
\int_C \int_{\R\times Y} \left|\mathbf 1_{\xi<u\left(t,x,y  \right)} - \mathbf 1_{\xi<u_{\delta}\left(t,x,y\right)}\right| \:K(\xi)\:dt\:dx\:dy\:d\xi= ||u-u_{\delta}||_{L^1(C\times Y)}.
\end{gather*}
Hence we have proved that for all $\delta>0$, for all compact set $C\subset [0,\infty)\times \R^N$,
$$
\lim_{\e\to 0}\left|\left| \ue(t,x) - u_{\delta} \left(t,x,\frac{x}{\e}  \right)\right|  \right|_{L^1(C)}= ||u-u_{\delta}||_{L^1(C\times Y)}.
$$

Statement \eqref{strongCVregul} then follows from standard convolution results.

\subsection{Application: proof of the convergence result for kinetic solutions}
In this subsection, we prove Theorem \ref{thm:stcvkinsol}; this result is in fact an easy consequence of the convergence result stated in Theorem \ref{thm:strg_cvg} for entropy solutions, and of the contraction principle for the limit system. Assume that $a_{N+1}\equiv 0$, and let $\ue$ be a kinetic solution of equation \eqref{eqdepart}, with an initial data $u_0(x,x/\e)$ such that $u_0\in L^1(\R^N, \mathcal C_{\text{per}}(Y))$ and
\be
\Si \frac{\p}{\p y_i}\left(a_i(y,\xi) \chi(\xi,u_0(x,y)) \right)=0\label{hyp:u_0wellprep}
\ee
in the sense of distributions.

For $n\in\N$, let $u_0^n:= \sgn(u_0) \inf(|u_0|,n)$. Then for all $n\in\N$, $u_0^n$ belongs to $L^\infty(\R^N\times Y)$ and
$$
u_0^n\to u_0\quad \text{as } n\to\infty\quad \text{in }L^1(\R^N, \mathcal C_{\text{per}}(Y)).
$$
Moreover, $\chi(\xi,u_0^n)=\chi(\xi,u_0) \mathbf 1_{\xi< n}$, and thus for all $n\in\N$, $u_0^n$ satisfies \eqref{hyp:u_0wellprep}.

For all $n,\e>0$, let $\ue_n\in\mathcal C([0,\infty), L^1(\R^N))\cap L^\infty([0,\infty)\times \R^N)$ be the unique entropy solution of equation \eqref{eqdepart} with initial data $u_0^n(x,x/\e)$. Then by the contraction principle for kinetic solutions of scalar conservation laws, we have
$$
\forall n\in\N,\quad ||\ue-\ue_n||_{L^\infty([0,\infty), L^1(\R^N))}\leq \left|\left| u_0\left( x,\frac{x}{\e} \right)-u_0^n\left( x,\frac{x}{\e} \right) \right| \right|_{ L^1(\R^N)}\leq ||u_0 - u_0^n||_{L^1(\R^N, \mathcal C_{\text{per}}(Y))}.
$$

On the other hand, for all $n\in\N$, let $\mathbf 1_{\xi<u_n}$ be the unique solution of the limit system with initial data $\mathbf 1_{\xi<u_0^n}$. By the contraction principle for solutions of the limit system (see inequality \eqref{in:contrac_zeta}), we have, for all integers $n,m\in\N$, for all $t\geq 0$,
\begin{eqnarray}
 \int_{\R^N\times Y}\left| u_n (t,x,y)- u_m(t,x,y) \right|\zeta(x)\:dx\:dy  &\leq& e^{Ct} \int_{\R^N\times Y}\left| u_0^m  (t,x,y)- u_0^n(t,x,y) \right|\zeta(x)\:dx\:dy \nonumber\\
 &\leq &e^{Ct} ||u_0^m - u_0^n||_{L^1(\R^N, \mathcal C_{\text{per}}(Y))},\label{Cauchy}
\end{eqnarray}
where $\zeta\in\mathcal C^\infty(\R^N)$ is a cut-off function satisfying the same hypotheses as in the previous subsections. Consequently, the sequence $(u_n)_{n\in\N}$ is a Cauchy sequence in $L^{\infty}_{\text{loc}}([0,\infty), L^1(\R^N\times Y))$; thus there exists a function $u\in L^{\infty}_{\text{loc}}([0,\infty), L^1(\R^N\times Y))$ such that $u_n$ converges towards $u$ as $n\to\infty$ in $L^{\infty}_{\text{loc}}([0,\infty), L^1(\R^N\times Y))$. Moreover, the limit $u$ is independent of the chosen sequence $u_0^n$ thanks to \eqref{Cauchy}: indeed, if $v_0^n$, $w^n_0$ are two approximating sequences giving rise to functions $v$ and $w$ respectively, we construct the sequence
$$
u_0^n=\left\{ \begin{array}{l}
            v_0^n\quad \text{if $n$ is even,}\\
            w_0^n\quad \text{if $n$ is odd.}\\
\end{array}\right.
$$
Then the sequence $u_0^n$ converges towards $u_0$, and thus the corresponding sequence $u_n$ converges towards $u$, while $u_{2n}$ converges towards $v$ and $u_{2n+1}$ towards $w$. By uniqueness of the limit, $u=v=w$.

On the other hand, since the sequence $\fe=\chi(\xi,\ue)$ is bounded in $L^\infty$, there exists a sequence $(\e_k)_{k\in\N}$ of positive numbers, $\e_k\to 0$, and a function $f\in L^{\infty}([0,\infty\times\R^N\times Y\times \R)$, such that
$$
\chi(\xi,u^{\e_k}(t,x))\stackrel{\text{2 sc.}}{\rightharpoonup}f(t,x,y,\xi).
$$
Now, for all $k,n\in\N$,
$$
||\chi(\xi,u^{\e_k})-\chi(\xi,u^{\e_k}_n)||_{L^\infty([0,\infty), L^1(\R^{N+1}))}\leq ||u_0 - u_0^n||_{L^1(\R^N, \mathcal C_{\text{per}}(Y))},
$$
and for all $n\in\N$, since $\chi(\xi,u)=\mathbf 1_{\xi<u}- \mathbf 1_{\xi<0}$, we have, as $k\to\infty$,
$$
\chi(\xi,u^{\e_k}_n)\stackrel{\text{2 sc.}}{\rightharpoonup}\chi(\xi,u_n).
$$
Let $\varphi \in\mathcal D_{\text{per}}([0,\infty)\times\R^N\times Y\times \R)$. By definition of two-scale convergence,
\begin{multline*}
\int_0^\infty\int_{\R^{N+1}}\left[\chi(\xi,u^{\e_k}(t,x))-\chi(\xi,u^{\e_k}_n(t,x))\right]\varphi \left( t,x,\frac{x}{\e_k},\xi \right)\:dt\:dx\:d\xi\to\\
\to \int_0^\infty\int_{\R^{N}\times Y\times \R}\left[f(t,x,y,\xi)-\chi(\xi,u_n(t,x,y))\right]\varphi \left( t,x,y,\xi \right)\:dt\:dx\:dy\:d\xi.
\end{multline*}
And for all $k\in\N$, the following inequality holds
\begin{multline*}
 \left|\int_0^\infty\int_{\R^{N+1}}\left[\chi(\xi,u^{\e_k}(t,x))-\chi(\xi,u^{\e_k}_n(t,x))\right]\varphi \left( t,x,\frac{x}{\e_k},\xi \right)\:dt\:dx\:d\xi\right| \leq\\
\leq ||\varphi||_{L^1([0,\infty), L^{\infty}(\R^N\times Y\times \R))}  ||u_0 - u_0^n||_{L^1(\R^N, \mathcal C_{\text{per}}(Y))}.
\end{multline*}

Passing to the limit as $k\to\infty$, we deduce that for all $n\in\N$, $\varphi \in\mathcal D_{\text{per}}([0,\infty)\times\R^N\times Y\times \R)$,
\begin{multline*}
 \left| \int_0^\infty\int_{\R^{N}\times Y\times \R}\left[f(t,x,y,\xi)-\chi(\xi,u_n(t,x,y))\right]\varphi \left( t,x,y,\xi \right)\:dt\:dx\:d\xi\right|\leq\\
\leq ||u_0 - u_0^n||_{L^1(\R^N, \mathcal C_{\text{per}}(Y))}||\varphi||_{L^1([0,\infty), L^{\infty}(\R^N\times Y\times \R)} .
\end{multline*}
Thus, we pass to the limit as $n\to\infty$ and we infer that $f=\chi(\xi,u(t,x,y))$ almost everywhere. Hence the limit is unique, and the whole sequence $\chi(\xi,\ue)$ converges (in the sense of two-scale convergence).

Eventually, let us pass to the limit as $n\to\infty$ in the limit evolution problem for $\chi(\xi,u_n)$. We set $f=\chi(\xi,u)$, and define the distribution
$$
\mathcal M:= \p_t f + a(y,\xi)\cdot \nabla_x f.
$$
Then $\mathcal M_n\rightharpoonup\mathcal M$ in the sense of distributions, and it is easily checked that inequality \eqref{hyp:Mdivfree} is preserved when passing to the (weak) limit. Thus $\mathcal M$ satisfies \eqref{hyp:Mdivfree}.
\vskip2mm

In the divergence-free case, the main difference between the $L^\infty$ and the $L^1$ setting, that is, Theorem~\ref{thm:strg_cvg} and Theorem \ref{thm:stcvkinsol}, lies in the fact that uniqueness for the limit system in the $L^1$ setting seems difficult to derive; indeed, the proof of uniqueness in the $L^\infty$ case uses several times the fact that the distribution $\mathcal M$ has compact support. In a $L^1$ setting, this assumption would have to be replaced by a hypothesis expressing that $\mathcal M$ vanishes as $|\xi|\to\infty$, in some sense. But it is unclear how to retrieve such a property from the hydrodynamic limit (see section \ref{sec:BGK}), for instance. The above argument only proves that uniqueness holds among $L^1$ solutions which are obtained as the limit of a sequence of $L^\infty$ solutions. Thus we have left open the correct notion of limit system in a weak $L^1$ setting, and the derivation of uniqueness therein.

Nonetheless, we wish to stress that the contraction principle in the $L^\infty$ setting is sufficient to ensure that the whole sequence $\chi(\xi,\ue)$ converges, even if uniqueness for the limit system fails.

\section{A relaxation model for the limit evolution problem\label{sec:BGK}}

In this section, we exhibit another way of finding solutions of the limit system in the divergence-free case. Indeed, the existence result in theorem  \ref{thm:lim_ev_eq} was proved by passing to the two-scale limit in \eqref{eq:depart_kf_gal}, and it may be interesting to have another way of constructing solutions, which does not involve a homogenization process.

Hence, we introduce a relaxation model of BGK type, in which we pass to the limit as the relaxation parameter goes to infinity. The drawback of this method is that the existence of solutions of the limit system is not a consequence of the construction. Indeed, we shall prove that if a solution of the limit system exists, then the family of solutions of the relaxation model strongly converge towards it in the hydrodynamic limit. Hence the proof is not self-contained, because the existence of a solution of the limit system is required in order to pass to the limit. Nevertheless, the final result may be useful in other applications.

In the whole section, when we refer to the limit system, we have in mind the modified equations introduced in Remark \ref{rem:lim_syst_divfree}. In the divergence-free case, it is also slightly more convenient to work with the function $\chi(\xi,u)$, rather than $\mathbf 1_{\xi<u}$. Hence a solution of the limit problem is a function $g$ satisfying
\begin{gather}
\Si \frac{\p}{\p y_i}\left(a_i(y,\xi) g \right)=0,\label{eq:cell_divfree}\\
\frac{\p g}{\p t}  + \Si a_i(y,\xi) \frac{\p g}{\p x_i}=\mathcal M,\\
\p_{\xi} g = \delta(\xi) - \nu(t,x,y,\xi), \quad \nu\geq 0,\end{gather}
and $\mathcal M$ is such that for all $\varphi\in\mathcal D([0,\infty)\times \R^N)$ such that $\varphi\geq 0$, the function $\mathcal M \ast_{t,x} \varphi$ belongs to $\mathcal C ([0,\infty)\times \R^N, L^2(Y\times \R))$, and
\be
\left\lbrace
\begin{array}{l}
\ds\int_{Y\times \R } \left(\mathcal M \ast_{t,x} \varphi\right)(t,x,\cdot) \:\psi\leq 0,\\
\ds\forall \psi\in L^{\infty}_{\text{loc}}(Y\times \R),\ \dv_y(a\psi)=0,\text{ and }\p_{\xi} \psi\geq 0.
\end{array}
\right.\label{hyp:M_divfree}
\ee

\subsection{A relaxation model\label{ssec:relax}}

The goal of this subsection is to introduce a system approaching \eqref{eq:cell_divfree}-\eqref{hyp:M_divfree}. With this aim in view, we define a relaxation model of BGK type, which takes into account the constraints the limit system, that is, equations \eqref{eq:cell_divfree}-\eqref{hyp:M_divfree}. Let
\begin{gather*}
M:=||u_0||_{L^\infty(Y\times \R)},\\
E:= \{f\in L²(Y\times \R) ,\  \supp f\subset Y\times [-M,M]\},\\
\mathbb K:=\left\lbrace  \varphi\in E, \ \dv_y (a(y,\xi) \varphi(y,\xi))=0 \quad \text{in }\mathcal D'.\right\rbrace ,\\
\mathcal K :=\mathbb K\cap \left\lbrace  \varphi\in E,\ \exists \nu\in M^1_\text{per}(Y\times \R), \ \nu\geq 0, \ \p_\xi \varphi=\delta(\xi) -\nu\right\rbrace .
\end{gather*}
Then $E$ endowed with the usual scalar product is a Hilbert space, and $\mathcal K$ is a closed convex set in $E$. Thus the projection $\mathcal P$ on $\mathcal K$ is well-defined.

The main result of this subsection is the following :
\begin{prop}Let $\lambda,T>0$ be arbitrary.
Set
$$
X_T:=\mathcal C([0,T], L^2(\R^N_x\times Y\times \R_\xi)).
$$

 Then there exists a unique solution $\fl\in X_T$ of the equation
\be\left\{
\begin{array}{l}
\p_t \fl + a(y,\xi)\cdot\nabla_x \fl + \lambda \fl = \lambda \mathcal P(\fl),\\
\fl(t=0,x,y,\xi)=\chi(\xi,u_0(x,y))
\end{array}\right.
\label{eq:BGK}
\ee
The function $\fl$ has the following properties :
\begin{enumerate}
 \item For almost every $t,x,y,\xi$,
\begin{gather*}
\fl(t,x,y,\xi)=0\text{ if }\xi\geq M,\\
\sgn(\xi)\fl(t,x,y,\xi)=|\fl(t,x,y,\xi)|\leq 1.
\end{gather*}
\item $L^2$ estimate: for all $\lambda>0$,
\be
||\fl||_{X_T}\leq ||u_0||_{L¹(\R^N\times Y)}.\label{est:fl}\ee

\item Strong continuity at time $t=0$: there exists a function $\omega:[0,\infty)\to[0,\infty)$, such that $\lim_{0^+}\omega=0$, and such that for all $\lambda>0$, $t\geq 0$,
\be
||\fl(t)-f_0||_{L^1(\R^N\times Y\times \R)}\leq \omega(t).\label{cont_time0_fl}
\ee

\item Fundamental inequality for $\mathcal M_{\lambda}:=\lambda (\mathcal P(\fl)-\fl)$: for all $g\in\mathcal K$, for almost every $(t,x)$,
\be
\int_{Y\times\R} \mathcal M_{\lambda}\:(\mathcal P(\fl)-g)\leq 0.\label{in:Ml}
\ee

\end{enumerate}

\label{prop:BGK}

\end{prop}

In equation \eqref{eq:BGK}, the projection $\mathcal P$ acts on the variables $y,\xi$ only; since $f$ is a function of $t,x,y,\xi$, $\mathcal P(f)$ should be understood as
$$
\mathcal P(f)(t,x,\cdot)=\mathcal P(f(t,x,\cdot)),
$$
and the above equality holds between functions in $L²(Y\times\R)$, almost everywhere in $t,x$.

\begin{proof}
{\it First step. Construction of $\fl$.} The existence and uniqueness of $\fl$ follows from a fixed point theorem in $X_T$. We define the application $\phi_T:X_Y\to X_T$ by $\phi_T(f)=g$, where $g$ is the solution of the linear equation
$$
\left\{
\begin{array}{l}
\p_t g + a(y,\xi)\cdot\nabla_x g + \lambda g = \lambda \mathcal P(f),\\
g(t=0,x,y,\xi)=\chi(\xi,u_0(x,y))
\end{array}\right.
$$
The existence and uniqueness of $g$ follows from well-known results on the theory of linear transport equations (recall that $a\in\mathcal C¹$). Moreover, if $f_1, f_2\in X_T$ and $g_i=\phi_T(f_i)$, $i=1,2$, then $g=g_1-g_2$ is a solution of
$$
\left\{
\begin{array}{l}
\p_t g + a(y,\xi)\cdot\nabla_x g + \lambda g = \lambda \left[\mathcal P(f_1)- \mathcal P(f_2)\right],\\
g(t=0,x,y,\xi)=0.
\end{array}\right.
$$
Multiplying the above equation by $g$, and integrating on $\R^N_x\times Y\times \R_\xi$, we obtain the estimate
$$
\frac 1 2 \frac{d}{dt}||g(t)||_{L²(\R^N\times Y\times\R)}^2 + \lambda ||g(t)||_{L²(\R^N\times Y\times\R)}^2\leq \lambda \int_{\R^N\times Y\times\R} \left[\mathcal P(f_1)- \mathcal P(f_2)\right] g.
$$
Recall that the projection $\mathcal P$ is Lipschitz continuous with Lipschitz constant 1. Thus
\begin{eqnarray*}
\int_{\R^N\times Y\times\R} \left[\mathcal P(f_1)- \mathcal P(f_2)\right] g&\leq& \frac 1 2 ||\mathcal P(f_1(t)) -\mathcal P(f_2(t))||_{L²(\R^N\times Y\times\R)}^2 + \frac 1 2 ||g(t)||_{L²(\R^N\times Y\times\R)}^2 \\
&\leq & \frac 1 2 ||(f_1-f_2)(t)||_{L²(\R^N\times Y\times\R)}^2 + \frac 1 2 ||g(t)||_{L²(\R^N\times Y\times\R)}^2.
\end{eqnarray*}
Eventually, we obtain
$$
\frac{d}{dt}||g(t)||_{L²(\R^N\times Y\times\R)}^2 + \lambda ||g(t)||_{L²(\R^N\times Y\times\R)}^2\leq \lambda  ||(f_1-f_2)(t)||_{L²(\R^N\times Y\times\R)}^2\leq \lambda ||f_1-f_2||_{X_T}^2.
$$
A straightforward application of Gronwall's lemma yields

$$
||g_1-g_2||_{X_T}\leq \sqrt{1-e^{-\lambda T}}\:||f_1-f_2||_{X_T}.
$$
Thus $\phi_T$ is a contractant application and has a unique fixed point in $X_T$, which we call $\fl$.

\vskip1mm

\noindent{\it Second step. $L^2$ estimate.} Multiplying \eqref{eq:BGK} by $\fl$ and integrating on $\R^N\times Y\times \R$, we infer
$$
\frac 1 2 \frac{d}{dt}||\fl(t)||_{L²(\R^N\times Y\times\R)}^2 + \lambda ||\fl(t)||_{L²(\R^N\times Y\times\R)}^2\leq \lambda \int_{\R^N\times Y\times \R}\mathcal P(\fl) \fl.
$$
Notice that $0\in\mathcal K$; thus the Lipschitz continuity of $\mathcal P$ entails that for almost every $t,x$
$$
||\mathcal P(\fl)(t,x)||_{E}\leq ||\fl(t,x)||_E.
$$
Hence, using the Cauchy-Schwartz inequality, we deduce that $t\mapsto ||\fl(t)||_{L²(\R^N\times Y\times\R)}$ is nonincreasing on $[0,T]$. The equality
$$
\int_{\R^N\times Y\times\R}|\chi(\xi,u_0(x,y))|²\:dx\:dy\:d\xi=\int_{\R^N\times Y\times\R}|\chi(\xi,u_0(x,y))|\:dx\:dy\:d\xi=\int_{\R^N\times Y}|u_0(x,y)|\:dx\:dy
$$
then yields the desired result.

\vskip1mm

\noindent{\it Third step. Compact support in $\xi$.} Let us prove now that $\fl(\cdot,\xi)=0$ if $|\xi| > M$: let $\varphi\in\mathcal D(\R)$ be an arbitrary test function such that $\varphi(\xi)=0$ when $|\xi|\leq M$. Then $\mathcal P(\fl) \varphi=0$ since $\mathcal P(\fl) \in\mathcal K$, and thus $\fl \varphi$ is a solution of
$$
\begin{array}{l}
\ds{\frac{\p}{\p t}\left( \fl \varphi\right) + a\cdot\nabla_x \left( \fl \varphi\right) } + \lambda\left( \fl \varphi\right)=0,\\
\left( \fl \varphi\right)(t=0,x,y,\xi)=0.
\end{array}
$$
Hence $(\fl \varphi)(t,x,y,\xi)=0$ for almost every $t,x,y,\xi$, and $\fl(\cdot,\xi)=0$ if $|\xi|>M$.

\vskip1mm

\noindent{\it Fourth step. Sign property.} We now prove the sign property, namely
$$
\sgn(\xi) \fl =|\fl|\leq 1\quad \text{a.e.}
$$
This relies on the following fact: if $g\in\mathcal K$, then $\sgn(\xi) g(y,\xi)\in[0,1]$ for almost every $y,\xi$. Indeed, $g(\cdot,\xi)=0$ if $\xi<-M$, and thus if $-M<\xi<0$,
$$
g(y,\xi) = -\int_{-M}^\xi \nu(y,\xi')\:d\xi'\leq 0.
$$
Hence $g(y,\cdot)$ is non positive and non increasing on $(-\infty, 0)$. Similarly, $g(y,\cdot)$ is non negative and non decreasing on $( 0,\infty)$. And if $\xi<0<\xi'$, then
$$
g(y,\xi') - g(y,\xi) = 1 - \int_{\xi}^{\xi'}\nu(y,w)\:dw \leq 1.
$$
Hence the sign property is true for functions in $\mathcal K$.

Multiplying \eqref{eq:BGK} by $\sgn(\xi)$, we are led to
$$
\frac{\p}{\p t}\left(\sgn(\xi) \fl\right) + a(y,\xi) \cdot\nabla_x\left(\sgn(\xi) \fl\right) + \lambda \left(\sgn(\xi) \fl\right) = \lambda \mathcal P(\fl)\in [0,\lambda].
$$
And at time $t=0$, $\sgn(\xi) \fl(t=0)=|\chi(\xi,u_0)|\in[0,1]$. Thus, using a maximum principle for this linear transport equation, we deduce that the sign property is satisfied for $\fl$.

\vskip1mm

\noindent{\it Fifth step. Uniform continuity at time $t=0$.} Let $\delta>0$ be arbitrary, and let $f_0^{\delta}:=f_0\ast_x \theta^\delta$, with $\theta^\delta$ a standard mollifier. Then $f_0^{\delta}(x)\in\mathcal K$ for all $x\in\R^N$, and thus $\fl-f_0^{\delta}$ is a solution of the equation
$$
\frac{\p}{\p t}\left ( \fl-f_0^{\delta}\right) + a(y,\xi)\cdot \nabla_x \left ( \fl-f_0^{\delta}\right)  + \lambda\left ( \fl-f_0^{\delta}\right) = \lambda\left ( \mathcal P(\fl)-\mathcal P(f_0^{\delta})\right) - a(y,\xi)\cdot \left ( f_0 \ast_x \nabla\theta^{\delta} \right).
$$
Multiply the above equation by $\left ( \fl-f_0^{\delta}\right)$ and integrate on $\R^N \times Y \times \R$. Using once more the Lipschitz continuity of the projection $\mathcal P$, we obtain
\begin{eqnarray*}
 \frac{1}{2}\frac{d}{dt}\left|\left|\fl - f_0^{\delta}\right|\right|_{L^2(\R^N \times Y \times \R)}^2&\leq &||a||_{L^{\infty}(Y\times(-M,M))}||\fl - f_0^{\delta}||_{L^2(\R^N \times Y \times \R)}||f_0||_{L^2(\R^N \times Y \times \R)}||\nabla\theta^{\delta}||_{L^1}\\
\frac{d}{dt}\left|\left|\fl - f_0^{\delta}\right|\right|_{L^2(\R^N \times Y \times \R)}&\leq & \frac{C}{\delta}.
\end{eqnarray*}
As a consequence, we obtain the following estimate, which holds for all $t>0$, $\lambda>0$ and $\delta>0$
$$
\left|\left|\fl(t)- f_0^{\delta}\right|\right|_{L^2(\R^N \times Y \times \R)}\leq \frac{C t}{\delta} + \left|\left|f_0- f_0^{\delta}\right|\right|_{L^2(\R^N \times Y \times \R)}.
$$
Hence the uniform continuity property is true, with
$$
\omega(t):=\inf_{\delta>0}\left (\frac{C t}{\delta} + 2\left|\left|f_0- f_0^{\delta}\right|\right|_{L^2(\R^N \times Y \times \R)} \right).
$$

\vskip1mm

\noindent{\it Sixth step. Inequality for $\mathcal M_{\lambda} $.} Inequality \eqref{in:Ml} is merely a particular case of the inequality
$$
\langle \mathcal P(f)-f,\mathcal P(f)-g\rangle_E\leq 0
$$
which holds for all $f\in E$, for all $g\in\mathcal K$.

\end{proof}

\subsection{The hydrodynamic limit\label{ssec:hydro}}

In this subsection , we prove the following result :
\begin{prop}
Let $(\fl)_{\lambda>0}$ be the family of solutions of the relaxation model \eqref{eq:BGK}, and let $f(t)=\chi(\xi,u)$ be the unique solution of the limit system \eqref{eq:cell_divfree}-\eqref{hyp:M_divfree} with initial data $\chi(\xi,u_0(x,y))$. Then as $\lambda\to \infty$,
$$
\fl \to f\quad \text{in } L^2((0,T)\times \R^N\times Y\times \R).
$$

\end{prop}

The above Proposition relies on an inequality of the type
$$
\frac{d}{dt}\int_{\R^N\times Y \times \R}|\fl - f|^2\leq r_{\lambda}(t),
$$
with $r_{\lambda}(t)\to 0$ as $\lambda\to \infty$. The calculations are very similar to those of the contraction principle in the previous section; the only difference lies in the fact that $\fl$ and $f$ are not solutions of the same equation.

Before tackling the proof itself, let us derive a few properties on the weak limit of the sequence $\fl$. Since the sequence $\fl$ is bounded in $X_T\subset L²((0,T)\times \R^N\times Y\times \R)$, we can extract a subsequence, which we relabel $\fl$, and find a function $g\in  L²((0,T)\times \R^N\times Y\times \R)$ such that $\fl$ weakly converges to $g$ in $L²$. Moreover, the sequence $\mathcal P(\fl)$ is bounded in $L²((0,T)\times \R^N\times Y \times \R)$, for all $T>0$. Hence, extracting a further subsequence if necessary, we can find a function $h\in L^2((0,T)\times \R^N\times Y \times \R) $ such that $\mathcal P(\fl)$ weakly converges towards $h$ as $\lambda\to\infty$. Notice that the convex set $\mathcal K$ is closed for the weak topology in $L²$. Consequently, $h(t,x)\in \mathcal K$ for almost every $t,x$. At last,
$$
\mathcal P(\fl) -\fl =\mathcal O \left(\frac{1}{\lambda}\right),
$$
where the $\mathcal O$ is meant in the sense of distributions. Hence, $g=h$, and in particular, we deduce that $g(t,x)\in \mathcal K$ for almost every $(t,x)$.

We are now ready to prove the contraction inequality; consider a mollifying sequence $\theta^{\delta}$ as in the previous section, and set
$f^{\delta}= f\ast_{t,x} \theta^{\delta}$, $\fl^{\delta'}= \fl\ast_{t,x} \theta^{\delta'}$. Then
\begin{gather*}
\p_t  f^{\delta} + a(y,\xi) \cdot \nabla_x f^{\delta} = \mathcal M^{\delta},\\
\p_t  \fl^{\delta'} + a(y,\xi) \cdot \nabla_x \fl^{\delta'} = \mathcal M_{\lambda}^{\delta'}.
\end{gather*}
Let us multiply the first equation by $ \sgn(\xi) - 2\fl^{\delta'} $, the second by $2(\fl^{\delta'}-f^{\delta})$, and add the two identities thus obtained; setting $F^{\delta, \delta'}_{\lambda}=\sgn(\xi) f^{\delta} + |\fl^{\delta'} |^2 - 2 f^{\delta} \fl^{\delta'}$, we have
$$
\p_t F^{\delta, \delta'}_{\lambda} + a(y,\xi)\cdot \nabla_x F^{\delta, \delta'}_{\lambda} = \mathcal M^{\delta}\left(\sgn(\xi) - 2\fl^{\delta'} \right) + 2 \mathcal M_{\lambda}^{\delta'}(\fl^{\delta'}-f^{\delta}).
$$
We integrate over $(0,t)\times\R^N\times Y \times \R$ and obtain
\begin{eqnarray*}
\int_{\R^N\times Y \times \R}F^{\delta, \delta'}_{\lambda}(t,x,y,\xi)\:dx\:dy\:d\xi &=& \int_0^t\int_{\R^N\times Y \times \R}\mathcal M^{\delta}\left(\sgn(\xi) - 2\fl^{\delta'} \right)\\
&&+ 2\int_0^t\int_{\R^N\times Y \times \R}\mathcal M_{\lambda}^{\delta'}(\fl^{\delta'}-f^{\delta})\\
&& + \int_{\R^N\times Y \times \R}F^{\delta, \delta'}_{\lambda}(t=0,x,y,\xi)\:dx\:dy\:d\xi .
\end{eqnarray*}
We now pass to the limit as $\delta'\to 0$, with all the other parameters fixed. Notice that
\begin{eqnarray*}
\lim_{\delta'\to 0 } \int_0^t\int_{\R^N\times Y \times \R}\mathcal M_{\lambda}^{\delta'}(\fl^{\delta'}-f^{\delta})&=&\int_0^t\int_{\R^N\times Y \times \R}\mathcal M_{\lambda}(\fl-f^{\delta})\\
&=&-\lambda \int_0^t\int_{\R^N\times Y \times \R}(\fl- \mathcal P(\fl))^2\\&& + \int_0^t\int_{\R^N\times Y \times \R}\mathcal M_{\lambda}(\mathcal P(\fl)-f^{\delta})\\
&\leq &0,
\end{eqnarray*}
since $f^{\delta}(t,x)\in\mathcal K$ for all $t,x$. The passage to the limit in $F^{\delta, \delta'}_{\lambda}(t=0)$ does not rise any difficulty because of the strong continuity of the functions $\fl$ at time $t=0$. Hence, we retrieve
\begin{eqnarray*}
&& \int_{\R^N\times Y \times \R} \left\{\left(| f^{\delta}(t)| -  |f^{\delta}(t)|^2 \right) + |f^{\delta}(t)- \fl(t)|^2  \right\}\\
&\leq &\int_0^t\int_{\R^N\times Y \times \R}\mathcal M^{\delta}\left(\sgn(\xi) - 2\fl\right) \\
&&+ \int_{\R^N\times Y \times \R} \left\{\left(| f^{\delta}(t=0)| -  |f^{\delta}(t=0)|^2 \right) + |f^{\delta}(t=0)- \chi(\xi,u_0)|^2  \right\},\end{eqnarray*}
and thus, integrating once again this inegality for $t\in[0,T]$,
\begin{eqnarray*}
 &&\int_0^T\int_{\R^N\times Y \times \R}\left\{\left(| f^{\delta}| -  |f^{\delta}|^2 \right) + |f^{\delta}(t)- \fl|^2  \right\}\\&\leq &\int_0^T dt\left[ \int_0^t \int_{\R^N\times Y \times \R}\mathcal M^{\delta}(s)\left(\sgn(\xi) - 2\fl(s)\right) \: ds\right]\\
&& + T \int_{\R^N\times Y \times \R} \left\{\left(| f^{\delta}(t=0)| -  |f^{\delta}(t=0)|^2 \right) + |f^{\delta}(t=0)- \chi(\xi,u_0)|^2  \right\}.
\end{eqnarray*}
We now pass to the limit as $\lambda\to\infty$, with $\delta >0$ fixed. Then
$$
\liminf_{\lambda \to\infty}||\fl- f^{\delta}||_{L^2((0,T)\times \R^N\times Y \times \R}^2\geq ||g- f^{\delta}||_{L^2((0,T)\times \R^N\times Y \times \R}^2,
$$
and
\begin{eqnarray*}
 &&\lim_{\lambda \to \infty }\int_0^T dt\left[ \int_0^t \int_{\R^N\times Y \times \R}\mathcal M^{\delta}(s)\left(\sgn(\xi) - 2\fl(s)\right) \: ds\right]\\
&=&\int_0^T dt\left[ \int_0^t \int_{\R^N\times Y \times \R}\mathcal M^{\delta}(s)\left(\sgn(\xi) - 2g(s)\right) \: ds\right]\leq 0.
\end{eqnarray*}
Thus, we obtain, for all $\delta>0$
$$
||g- f^{\delta}||_{L^2((0,T)\times \R^N\times Y \times \R}^2\leq  T \int_{\R^N\times Y \times \R} \left\{\left(| f^{\delta}(t=0)| -  |f^{\delta}(t=0)|^2 \right) + |f^{\delta}(t=0)- \chi(\xi,u_0)|^2  \right\}.
$$
We have already proved in the previous section that the family $f^{\delta}(t=0)$ strongly converges towards $\chi(\xi,u_0)$ as $\delta $ vanishes, due to the continuity assumption at time $t=0$. Hence, we obtain in the limit
$$
||g- f||_{L^2((0,T)\times \R^N\times Y \times \R}^2\leq 0,
$$
and consequently, $g=f$. Hence the result is proved.

\section{The separate case : identification of the limit problem\label{sec:separated}}

This section is devoted to the proof of Proposition
\ref{thm:homogenized_sep}. Thus we focus on the limit system in the case where the flux $A$ can be written as
$$
A(y,\xi)=a_0(y) g(\xi),\quad \text{with }\dv_y a_0=0.
$$
The interest of this case lies in the special structure of the limit system; indeed, we shall prove that the function $u$, which is the two-scale limit of the sequence $\ue$, is the solution of the scalar conservation law \eqref{eq:homogenized_sep}. In view of Theorem \ref{thm:lim_ev_eq},
we wish to emphasize that Proposition \ref{thm:homogenized_sep}
implies in particular that the entropy solution of
\eqref{eq:homogenized_sep} satisfies the constraint equation
$$
\dv_y \left(a_0(y) g(u(t,x;y)) \right)
$$
for almost every $t>0,x\in\R^N$; this fact is not completely
obvious when $g\neq \Id$. We will prove in the sequel that
$u(t,x)$ actually belongs to the constraint space $\K_0$ for a.e.
$t,x$.

Before tackling the proof of the theorem, let us mention that the
limit problem \eqref{eq:homogenized_sep} is not the one which is expected from a vanishing viscosity approach. Precisely, for
any given $\delta>0$, let $\ue_{\delta}$ be the
solution of
$$
\p_t \ue_{\delta} + \dv_x A\left(\frac x {\e},\ue_{\delta}
\right)-\e\delta\Delta_x\ue_{\delta}=0,
$$
with the initial data $\ue_{\delta}(t=0,x)=u_0\left(x,x/\e
\right)$. Then for all $\e>0$, $\ue_{\delta}\to \ue$ as $\delta\to
0$; moreover, the behavior of $\ue_{\delta}$ as $\e\to 0$ is known
for each $\delta>0$ (see \cite{homogpara,initiallayer}). In the divergence-free case, for all
$\delta>0$,
$$
\lim_{\e\to 0}\ue_{\delta}=\bu(t,x)\quad \text{in } L^1_{\text{loc}},
$$
where $\bu$ is the entropy solution of
$$
\p_t \bu + \dv_x (\mean{a} g(\bu))=0,
$$
with initial data $\bu(t=0,x)=\mean{u_0(x,\cdot)}$. Hence, it could be expected that the limits $\e\to 0$ and $\delta\to 0$ can be
commuted, that is
$$
\lim_{\e\to 0}\lim_{\delta\to 0}\ue_{\delta}=\lim_{\delta\to
0}\lim_{\e\to 0}\ue_{\delta},
$$
which would entail
$$
\lim_{\e\to 0}\ue =\bu.
$$

In general, this equality is false, even in a weak sense: a
generic counter-example is the one of shear flows (see for
instance the calculations in \cite{wetransport}). In that case, we
have $N=2$ and $A(y,\xi)=(a_1(y_2)\xi,0)$, and the equation
\eqref{eq:homogenized_sep} becomes
$$
\p_t u + a_1(y_2)\p_{x_1}u=0,
$$
with the initial condition $u(t=0,x,y)=u_0(x_1,x_2,y_2)$. It is
then easily checked that in general, the average of $u$ over $Y$
is not the solution of the transport equation
$$
\p_t \bu + \mean{a_1}\p_{x_1}\bu=0.
$$

We now turn to the proof of Proposition \ref{thm:homogenized_sep}. In view of Theorem \ref{thm:lim_ev_eq}, it is sufficient to prove that the entropy solution
of \eqref{eq:homogenized_sep} belongs to $\K_0$ for a.e. $t,x$, or
in other words, that $\K_0$ is invariant by the semi-group associated to equation
\eqref{eq:homogenized_sep}. We prove this result in the slightly more general context of kinetic solutions. The core of the proof lies in the following

\begin{prop}
Let $u_0\in L^1(\R^N,L^{\infty}(Y))$ such that
$u_0(x,\cdot)\in\K_0$ for almost every $x\in\R^N$.

Let $v=v(t,x;y)\in \mathcal C([0,\infty); L^1(\R^N\times Y))$ be
the kinetic solution of
$$
\left\{ \begin{array}{l} \p_t v(t,x;y) + \dv_x\left(\ta(y)
g(v(t,x;y)) \right)=0,\quad t>0,\ x\in\R^N,\ y\in Y,\\
v(t=0,x;y)=u_0(x,y),
\end{array}\right.
$$
i.e. $f^1(t,x,y,\xi):=\chi(\xi,v(t,x;y))$ is a solution in the
sense of distributions of \be \left\{ \begin{array}{l} \p_t f^1 +
\ta(y)\cdot \nabla_x f^1 g'(\xi)
=\p_\xi m,\quad t>0,\ x\in\R^N,\ y\in Y,\ \xi\in\R,\\
f^1(t=0,x,y,\xi)=\chi(\xi,u_0(x,y)),
\end{array}\right.
\label{eq:f1}\ee and $m$ is a non-negative measure on
$[0,\infty)\times\R^N\times Y\times \R$.

Then for a.e. $t>0, x\in\R^N, \ u(t,x)\in \K_0$.

\label{prop:hydro}

\end{prop}

\begin{proof}

First, let us recall (see \cite{BP,PT}) that for all
$T>0$,
$$
f^1=\lim_{\lambda\to\infty}f_{\lambda}\quad\text{in }\mathcal
C([0,T]; L^1(\R^N\times Y\times\R)),
$$
where $f_{\lambda}=f_{\lambda}(t,x,y,\xi)$ ($\lambda>0$) is the
unique solution of the system \be \left\{
\begin{array}{l}
\p_t \fl + \ta(y)\cdot\nabla_x \fl\: g'(\xi) + \lambda\fl
=\lambda\chi(\xi,\ul),\\
\ul(t,x,y)=\int_{\R}\fl(t,x,y,\xi)\:d\xi,\\
\fl(t=0)=\chi(\xi,u_0).
\end{array}
\right. \label{eq:hydrolim}\ee

\noindent Moreover, for every $\lambda>0$, $\ul$ is the unique fixed point
of the contractant application
$$
\phi_{\lambda}:\begin{array}{rcl}
            \mathcal C((0,T); L^1(\R^N\times Y))&\to &\mathcal C((0,T); L^1(\R^N\times Y))\\
            u_1 &\mapsto & u_2
            \end{array}
$$
where $u_2=\int_{\xi} f$ and $f$ is the solution of \be
\begin{array}{l}
 \ds \p_t f +
\ta(y)\cdot\nabla_x f \:g'(\xi) + \lambda f =\lambda\chi(\xi,u_1),\\
f(t=0)=\chi(\xi,u_0).
\end{array}
\label{eq:evolhydro}\ee
Thus, the whole point is to prove that the space
$$
\{u\in \mathcal C([0,T]; L^1(\R^N\times Y)); u(t,x)\in\K_0 \ \text{a.e} \}
$$
is invariant by the application $\phi_{\lambda}$.

First, let us stress that for all $u\in L^{1}(Y)$,
\be
u\in\K_0\iff \dv_y (a(y)\chi(\xi,u))=0 \text{ in }\mathcal D'(Y\times\R).\label{equiv:uK0}
\ee
Indeed, if $u\in\K_0$, then for all $\delta >0$, set $u_{\delta}= u\ast \theta^{\delta}$, with $\theta^{\delta}$ a standard mollifier. The function $u_{\delta}$ is a solution of
$$
\dv_y (a_0 u_{\delta})= r_{\delta},
$$
and the remainder $r_{\delta}$ vanishes strongly in $L^1(Y)$ (see the calculations in the previous sections). Since the function $u_{\delta}$ is smooth, if $G\in \mathcal C^1(\R^N)$, we have
$$
\dv_y (a_0 G(u_{\delta}))= G'(u_{\delta})r_{\delta}.
$$
Passing to the limit as $\delta $ vanishes, we infer $\dv_y (a_0 G(u))=0$ for all $G \in\mathcal C^1(\R^N)$. At last, taking a sequence of smooth functions approaching $\chi(\xi,u)$, we deduce that $\dv_y(a_0\chi(\xi,u))=0$ in $\mathcal D'_{\text{per}}(Y\times \R)$. Conversely, assume that $\dv_y(a_0\chi(\xi,u))=0$; then integrating
this equation with respect to $\xi$ yields $u\in\K_0$. Hence \eqref{equiv:uK0} is proved.

Now, let $u_1\in C([0,T]; L^1(\R^N\times Y))$ such that
$u_1(t,x)\in\K_0$ a.e. Then $\dv(a_0\chi(\xi,u_1)=0)$. Let $f$ be the
solution of \eqref{eq:evolhydro}; since $\ta\in \K_0$, the distribution $\dv_y(a_0f)$
satisfies the transport equation
$$
\p_t \left( \dv(a_0f)\right) + g'(\xi)\ta(y)\cdot \nabla_x \left(
\dv(a_0f)\right) + \lambda\dv(a_0f)=0,
$$
and $\dv(a_0f)(t=0)=0$ because $u_0(x)\in\K_0$ a.e. Hence
$\dv_y(a_0f)=0$; integrating this equation with respect to $\xi$ gives
$u_2\in\K_0$ a.e.

Consequently, $\ul(t,x;\cdot)\in\K_0$ a.e. Passing to the limit,
we deduce that $v(t,x;\cdot)\in\K_0$ a.e.

\end{proof}

Let us now re-write equation \eqref{eq:f1}: setting $b(y)=a_0(y)-\ta(y)$, we have
$$
\p_t f^1 + a_0(y) \nabla_x f^1 g'(\xi) = \p_{\xi} m - b(y) \nabla_x f^1 g'(\xi)=: \mathcal M_1.
$$
If $u_0\in L^\infty(\R^N)$, then $v\in L^\infty([0,\infty)\times \R^N\times Y)$, and it is easily checked that $f^1$ and $\mathcal M_1$ satisfy the compact support assumptions. According to the above Proposition, $f^1$ also satisfies \eqref{eq:cell_divfree}, and thanks to the structure of the right-hand side, the distribution $\mathcal M_1$ satisfies \eqref{hyp:M_divfree}. Thus $f^1$ is the unique solution of the limit system, and Proposition \ref{thm:homogenized_sep} is proved.


\section{Further remarks on the notion of limit system\label{sec:furtherremarks}}

Here, we have gathered, by way of conclusion, a few remarks around the limit evolution system introduced in definition \ref{def:lim_ev_eq}. The main idea behind this section is that the limit system is not unique (although its solution always is), and thus several other relevant equations can be written instead of \eqref{eq:lim_evol_eq}. Unfortunately, there does not seem to be any rule which would allow to decide between two limit systems.

Let us illustrate these words by a first series of examples : assume that the flux is divergence free, and let
$$
\K:=\{f\in L^1_{\text{loc}}(Y\times \R),\ \Si \p_{y_i}(a_i f)=0\quad\text{in
}\mathcal D'\}.
$$
We denote by $P$ the projection on $\K$ in
$L^1_{\text{loc}}(Y\times\R)$. Precisely, consider the dynamical
system $X(t,y;\xi)$ defined by
$$
\left\{\begin{array}{l}
        \dot{X}(t,y;\xi)=a(X(t,y;\xi),\xi),t>0\\
        X(t=0,y;\xi)=y.
        \end{array}\right.
$$

Then for all $\xi\in\R$, the Lebesgue measure on $Y$ is invariant
by the semi-group $X(t;\xi)$ because of the hypothesis $\dv_y
a(y,\xi)=0$. Hence by the ergodic theorem, for all $f\in
L^1_{\text{loc}}(Y\times \R)$, there exists a function in
$L^1_{\text{loc}}(Y\times \R)$, denoted by $P(f)(y,\xi)$, such
that
$$
P(f)(y,\xi)=\lim_{T\to\infty}\frac 1 T\int_0^T
f(X(t,y;\xi),\xi)\:dt,
$$
and the limit holds a.e. in $y,\xi$ and in $Y\times (-R,R)$ for
all $R>0$.

Set $\tilde a : = P(a)$. Then if $f$ is a solution of the limit system, $f$ also satisfies
$$
\p_t f + \tilde a(y,\xi)\cdot \nabla_x f = \tilde{\mathcal M}
$$
and $f$, $\tilde{\mathcal M}$ satisfy \eqref{eq:cell_kin} and \eqref{hyp:dfdxi} - \eqref{hyp:M}. Indeed,
$$
\tilde{\mathcal M}= \mathcal M + \left[\tilde a(y,\xi) - a(y,\xi)\right]\cdot \nabla_x f
$$
and the term $\left[\tilde a(y,\xi) - a(y,\xi)\right]\cdot \nabla_x (f\ast_x \varphi)(t,x,y,\xi)$ belongs to $\K^\bot$ for all $t,x$. Of course, uniqueness holds for this limit system (the proof is exactly the same as the one in section \ref{sec:contraction}), and thus this constitutes as legitimate a limit system as the one in definition \ref{def:lim_ev_eq}. In fact, in the separate case, Proposition \ref{thm:homogenized_sep} indicates that the above system seems to be the relevant one, rather than the one in definition \ref{def:lim_ev_eq}. Notice that the distribution $\tilde{\mathcal M}$ satisfies the additional property
$$
\tilde{\mathcal M} \ast_{t,x} \phi(t,x) \in \K^\bot\quad \forall t,x.
$$

Let us now go a little further: let $\theta \in \mathcal C^1 (\R)$ such that $0\leq \theta \leq 1$, and let
$$a_{\theta}(y,\xi) =  \theta(\xi) a(y,\xi) + (1-\theta(\xi)) \tilde a(y,\xi).
$$
Then $f$ is a solution of
$$
\p_t f + \tilde a_{\theta}(y,\xi)\cdot \nabla_x f = \mathcal M_{\theta},
$$
for some distribution $\mathcal M_{\theta}$ satisfying \eqref{hyp:M}. Thus this still constitutes a limit system which has the same structure as the one of definition \ref{def:lim_ev_eq}. Hence the limit system is highly non unique, and it must be seen as a way of identifying the two-scale limit of the sequence $\fe$, rather than as a kinetic formulation of a given conservation law, for instance. We wish to emphasize that if the flux $A$ is not ``separated'', that is, if the hypotheses of Proposition \ref{thm:homogenized_sep} are not satisfied, then in general, the function $u$ such that $f=\mathbf 1_{\xi<u}$ is a solution of the limit system, is different from the solution $v=v(t,x,y)$ of the scalar conservation law
$$
\p_t v + \dv_x \tilde{A} (y,v)=0,
$$
where the flux $\tilde{A} $ is such that $\p_{\xi }\tilde{A_i} (y,\xi)= \tilde a_i(y,\xi)$. Indeed, the function $v$ above is not a solution of the cell problem in general, even if $v(t=0)$ is. In other words, the set $\K$ is not invariant by the evolution equation
$$
\p_t g +\Si \tilde a_i(y,\xi)\p_{x_i} g = \p_{\xi} m,
$$
where $m$ is a non-negative measure and $g=\mathbf 1_{\xi<v}$.

\vskip2mm

Let us now assume that the flux $A$ is not divergence free. Then there are cases where yet another notion of limit problem can be given: assume that there exists real numbers $p_1<p_2$, and a family $\{v(\cdot,p)\}_{p_1\leq p\leq p_2}$, which satisfies the following properties:
\begin{enumerate}
\item The function $(y,p)\mapsto v(y,p)$ belongs to $L^{\infty}(Y\times [p_1,p_2])$;
 \item For all $p\in[p_1,p_2]$, $v(\cdot,p)$ is an entropy solution of the cell problem; in other words, there exists a nonnegative measure $m(y,\xi;p)$ such that $f(y,\xi;p)=\mathbf 1_{\xi<v(y,p)}$ is a solution of
$$
\Si \frac{\p}{\p y_i}\left(a_i(y,\xi) f \right) + \frac{\p}{\p \xi}\left(a_{N+1}(y,\xi) f \right)=\frac{\p}{\p\xi} m;
$$
\item For all $p\in[p_1,p_2]$, $\mean{v(\cdot,p)}_Y=0$;

\item The distribution $\p_p v$ is a nonnegative function in $L¹(Y\times [p_1,p_2])$; this implies in particular that for all couples $(p,p')\in[p_1,p_2]^2$ such that $p\geq p'$, for almost every $y\in Y$,
$$
v(y,p)\geq v(y,p').
$$

\end{enumerate}

Under these conditions, one can construct a kinetic formulation for equation \eqref{eqdepart}, based on the family $v(x/\e,p)$ of stationary solutions of \eqref{eqdepart}, rather than on the family of Kruzkov's inequalities. This kind of construction was achieved in \cite{kinformpara} in a parabolic setting, following an idea developed by Emmanuel Audusse and Benoît Perthame in \cite{AudussePerthame}; these authors define a new notion of entropy solutions for a heterogeneous conservation law in dimension one, based on the comparison with a family of stationary solutions. Let us explain briefly how the kinetic formulation for entropy solutions of \eqref{eqdepart} is derived: let $\ue$ be an entropy solution of \eqref{eqdepart}. Define the distribution $\me\in\mathcal D'((0,\infty)\times \R^N\times (p_1,p_2))$ by
\be
\me(t,x,p):=-\left\{\frac{\p}{\p t}\left(\ue - v\left( \frac{x}{\e},p\right)\right)_+ +  \frac{\p}{\p y_i}\left[\mathbf 1_{ v\left( \frac{x}{\e},p\right)<\ue} \left(A_i\left(\frac{x}{\e},\ue\right) -A_i\left(\frac{x}{\e},v\left( \frac{x}{\e},p\right) \right) \right)\right] \right\}.
\label{def:me_kf}\ee
Then according to the comparison principle (which was known by Kruzkhov, see \cite{kruzkhov1,kruzkhov2}), $\me$ is a nonnegative measure on $(0,\infty)\times \R^N \times [p_1,p_2] $. Now, set
$$
\fe(t,x,p):= \mathbf 1_{ v\left( \frac{x}{\e},p\right)<\ue(t,x)}\in L^\infty([0,\infty)\times \R^N \times [p_1,p_2]).
$$
Thanks to the regularity assumptions on the family $v(\cdot,p)$, we can differentiate equality \eqref{def:me_kf} (which is meant in the sense of distributions) with respect to $p$, and we are led to
\be
\frac{\p}{\p t}\left(\fe v_p \left( \frac{x}{\e},p\right)\right) + \frac{\p}{\p x_i}\left(\fe v_p \left( \frac{x}{\e},p\right) a_i\left(\frac{x}{\e}, v \left( \frac{x}{\e},p\right)\right)\right)=\frac{\p \me}{\p p}.\label{eq:kf_heterogeneous}
\ee
This equation is in fact the appropriate kinetic formulation in the heterogeneous case; its main advantage on the equation \eqref{eq:depart_kf_gal} is the absence of the highly oscillating term
$$
\frac{1}{\e}\p_{\xi}\left[ a_{N+1}\left(\frac{x}{\e},\xi  \right) \mathbf1_{\xi<\ue}\right].
$$
Notice that for all $p\in[p_1,p_2]$,
\be
\dv_y\left(\frac{\p v(y,p)}{\p p} a(y,v(y,p))\right)=0\quad \text{in }\mathcal D'_{\text{per}}(Y ).\label{eq:dvdp}
\ee
This equation is derived by differentiating equation
$$
\dv_y A(y,v(y,p))=0
$$
with respect to $p$. Thus, if we set
$$
\tilde{a}(y,p):=\frac{\p v(y,p)}{\p p} a(y,v(y,p)),
$$
the vector field $\tilde{a}\in L^1(Y\times [p_1,p_2])$ is divergence-free, and the same kind of limit system as in the divergence free cas can be made. Of course, the interest of such a construction lies in the simplicity of the structure of the limit system in the divergence free case.

\begin{defi}

Let $f\in L^{\infty}([0,\infty), L^1(\R^N\times Y \times\R))$, $u_0\in L^1 \cap L^{\infty}(\R^N\times Y)$. We say that $f$ is a {\it generalized kinetic solution of the limit problem} associated with the family $v(\cdot,p)$ if there exists a distribution $\mathcal M\in \mathcal D'_{\text{per}}([0,\infty)\times \R^N \times Y \times \R)$ such that $f$ and $\mathcal M$ satisfy the following properties:

\begin{enumerate}
 \item Compact support in $p$: there exists $(p_1',p_2')\in[p_1,p_2]²$ such that $p_1<p_1'\leq p_2'<p_2$, and
\begin{gather*}
\supp \mathcal M \subset [0,\infty)\times \R^N \times Y \times [p_1',p_2'];\\
f(t,x,y,p)=1\text{ if }p_1<p<p_1', \quad f(t,x,y,p)=0\text{ if }p_2'<p<p_2.
\end{gather*}

\item Microscopic equation for $f$: $f$ is a solution in the sense of distributions on $Y\times (p_1,p_2)$ of
\be
\dv_y (\tilde{a}(y,p) f(t,x,y,p))=0.\label{eq:cell_kin2}\ee
\item Evolution equation: the couple $(f,\mathcal M)$ is a solution in the sense of distributions on $[0,\infty)\times \R^N\times Y\times (p_1,p_2)$ of
\be
\left\lbrace
    \begin{array}{l}
        \p_t (v_p(y,p) f) + \tilde a(y,p)\cdot \nabla_x f = \mathcal M,\\
        f(t=0,x,y,p)=\mathbf 1_{v(y,p)<u_0(x,y)}=:f_0(x,y,p);
\end{array}
\right.
\label{eq:lim_evol_eq2}
\ee
In other words, for any test function $\phi\in\mathcal D_{\text{per}}([0,\infty)\times \R^N \times Y \times (p_1,p_2))$,
\begin{multline*}
\int_0^\infty \int_{\R^N\times Y\times\R}f(t,x,y,p)v_p(y,p)\left\lbrace \p_t \phi(t,x,y,p)+ a(y,v(y,p))\cdot \nabla_x \phi (t,x,y,p)\right\rbrace \:dt\:dx\:dy\:d\xi =\\
=-\left\langle \phi,\mathcal M\right\rangle_{\mathcal D,\mathcal D'} -  \int_{\R^N\times Y\times\R}\mathbf 1_{v(y,p)<u_0(x,y)}v_p(y,p)\phi(t=0,x,y,p)\:dx\:dy\:d\xi.
\end{multline*}
\item Conditions on $f$: there exists a nonnegative measure $\nu\in M^1_{\text{per}}([0,\infty)\times \R^N \times Y \times \R)$ such that
\begin{gather}
\p_{p} f = -\nu,\label{hyp:dfdp}\\
0\leq f(t,x,y,\xi)\leq 1\quad \text{a.e.},\\
\frac{1}{\tau}\int_0^\tau\left|\left|f(s)-f_0 \right|  \right|_{L^2(\R^N\times Y\times (p_1,p_2)}\:ds\underset{\tau\to 0}{\longrightarrow}0\label{hyp:cont_time02}.
\end{gather}

\item Condition on $\mathcal M$: for all $\varphi\in\mathcal D([0,\infty)\times \R^N$ such that $\varphi\leq 0$, the function $\mathcal M \ast_{t,x} \varphi$ belongs to $\mathcal C ([0,\infty)\times \R^N, L^1(Y\times \R))$, and
\be
\left\lbrace
\begin{array}{l}
\int_{Y\times \R } \left(\mathcal M \ast_{t,x} \varphi\right)(t,x,\cdot) \:\psi\leq 0,\\
\forall \psi\in L^{\infty}_{\text{loc}}(Y\times \R),\ \dv_y(\tilde a\psi)=0,\text{ and }\p_{\xi} \psi\geq 0.
\end{array}
\right.
\label{hyp:M2}
\ee

\end{enumerate}

\label{def:lim_ev_eq2}
\end{defi}

We now state without proof a result analogue to Theorems \ref{thm:lim_ev_eq}, \ref{thm:strg_cvg} :

\begin{prop}
Let $A\in W^{2,\infty}_{\text{per,loc}}(Y\times \R)$. Assume that $a\in \mathcal C^1_{\text{per}}(Y\times \R)$ and that $\tta\in W^{1,1}(Y\times (p_1,p_2))$. Let $u_0\in  L^\infty(\R^N\times Y)\cap L^1_\text{loc}(\R^N,\mathcal C_{\text{per}}(Y))$ such that $u_0(x,\cdot)$ is an entropy solution of the cell problem for almost every $x\in\R^N$. Assume furthermore that there exists $p_1'<p_2'$ in $(p_1,p_2)^2$ such that
$$
v(y,p_1')\leq u_0(x,y)\leq v(y,p_2'),
$$
and let
$$
f_0(x,y,p):= \mathbf 1_{v(y,p)<u_0(x,y)}
$$

Then the following results hold :
\begin{enumerate}
 \item There exists a unique generalized kinetic solution $f$ of the limit problem associated with the family $(v(\cdot,p))_{p_1\leq p\leq p_2}$ with initial data $f_0$. Moreover, there exists a function $u\in L^\infty([0,\infty)\times \R^N\times Y)$ such that $$f(t,x,y,p)=\mathbf 1_{v(y,p)<u(t,x,y)}\quad \text{a.e.}$$

\item Let $\ue\in L^\infty([0,\infty)\times \R^N)$ be the entropy solution of \eqref{eqdepart} with initial data $u_0\left(x,x/\e \right)$. Let $f(t,x,y,p)=\mathbf 1_{v(y,p)<u(t,x,y)}$ be the unique solution of the limit problem. Then for all regularization kernels
$\varphi^{\delta}$ of the form
$$\varphi^{\delta}(x)=\frac 1{\delta^N} \varphi\left(\frac x{\delta} \right), \quad x\in\R^N,$$
with $\varphi\in\mathcal D(\R^N)$, $\int \varphi=1$, $0\leq
\varphi\leq 1$, we have, for all compact
$K\subset[0,\infty)\times\R^N$, \be \lim_{\delta\to 0} \lim_{\e\to
0}\left|\left |\ue(t,x) - u\ast_x\varphi^{\delta}\left(t,x,\frac
x\e \right)\right| \right|_{L^1(K)}=0.\ee

\end{enumerate}

\end{prop}

Hence a whole variety of limit systems can be given, depending on the choice of the family of solutions of the cell problem. However, it is not obvious that any given system is ``better'' than another one. But the important result, as far as homogenization is concerned, is that all systems have a unique solution.

\bibliography{../../articles,../../books}

\noindent CEREMADE-UMR 7534\\ Université Paris-Dauphine\\
Place du maréchal de Lattre de Tassigny\\ 75775 Paris Cedex 16, FRANCE\\
email: {\texttt{dalibard@ceremade.dauphine.fr}}

\end{document}